\documentclass[11pt,oneside]{article}
\usepackage{amsfonts}
\usepackage{graphics}
\usepackage{amsmath}
\usepackage{latexsym}
\usepackage{amssymb}
\usepackage{txfonts}
\usepackage{mathrsfs,tikz}
\usepackage{ifthen}
\usepackage{array}
\usepackage{cite}
\newtheorem{thm}{Theorem}[section]
\newtheorem{lem}{Lemma}[section]
\newtheorem{cor}{Corollary}[section]
\newtheorem{rmk}{Remark}[section]

\newtheorem{con}{Conjecture}[section]

\newtheorem{defn}{Definition}[section]

\newcommand{\NI}{\noindent}
\usepackage{authblk}
\title {On Hamilton Cycle Decompositions of Tensor Products of Graphs}
\author[1]{P. Paulraja\thanks{ppraja56@gmail.com}}
\author[2]{S. Sampath Kumar\thanks{ssamkumar.2008@gmail.com}}
\affil[1]{\small Department of Mathematics, Annamalai University\\Annamalainagar-608002, India}
\affil[2]{\small Department of Mathematics, SSN College of Engineering\\Kalavakkam-603110, India}

\begin{document}
 \date{}
\maketitle
\begin{abstract}
A Hamiltonian decomposition of $G$ is a partition of its edge set into disjoint Hamilton cycles. Manikandan and Paulraja conjectured that if $G$ and $H$ are Hamilton cycle decomposable circulant graphs with at least one of them is nonbipartite, then their tensor product is Hamilton cycle decomposable. In this paper, we have proved that, if $G$ is a Hamilton cycle decomposable circulant graph with certain properties and $H$ is a Hamilton cycle decomposable multigraph, then their tensor product is Hamilton cycle decomposable. In particular, tensor products of certain sparse Hamilton cycle decomposable circulant graphs are Hamilton cycle decomposable.
\end{abstract}
\section{\bf Introduction.}
All graphs considered here are simple and finite, unless otherwise stated. Let $C_k\,(\text{respectively }P_{k}),$ denote the cycle (respectively path) on $k$ vertices. The complete graph on $n$ vertices is denoted by $K_n$ and its complement is denoted by $\overline{K}_n.$ For two graphs $G$ and $H$ their {\it tensor product}, denoted by $G\times H,$ has vertex set $V(G)\times V(H)$ in which two vertices  $(g_1,h_1)\text { and }(g_2,h_2)$ are adjacent whenever $g_1g_2\in E(G)\text{ and }h_1h_2\in E(H).$ The {\it wreath product} of the graphs $G$ and $H,$ denoted by $G\circ H,$ has vertex set $V(G)\times V(H)$ in which $(g_1,\,h_1)(g_2,\,h_2)$ is an edge whenever $g_1g_2$ is an edge in $G,$ or $g_1=g_2$ and $h_1h_2$ is an edge in $H.$ Similarly, the \emph{cartesian product} of the graphs $G$ and $H,$ denoted by $G\Box H,$ has vertex set $V(G)\times V(H)$ in which $(g_1,\,h_1)(g_2,\,h_2)$ is an edge whenever $g_1=g_2$ and $h_1h_2$ is an edge in $H,$ or $h_1=h_2$ and $g_1g_2$ is an edge in $G.$ The subgraph induced by $S\subseteq V(G)$ is denoted by $\langle S\rangle.$ Similarly, the subgraph induced by $E^\prime\subseteq E(G)$ is denoted by $\langle E^\prime\rangle.$ For a graph $G,$ $G(\lambda)$ is the graph obtained from $G$ by replacing each edge of $G$ by $\lambda$ parallel edges. For a graph $G,$ $G^*$ is the symmetric digraph of $G.$

For two loopless multigraphs $G(\lambda)$ and $H(\mu),$ their tensor product, denoted by $G(\lambda)\times H(\mu),$ has the vertex set $V(G)\times V(H)$ and its edge set is described as follows: if $e=g_1g_2$ is an edge of multiplicity $\lambda$ in $G(\lambda)$ and $f=h_1h_2$ is an edge of multiplicity $\mu$ in $H(\mu),$ then corresponding to these edges there are edges $(g_1,h_1)(g_2,h_2)$ and $(g_1,h_2)(g_2,h_1)$ each of multiplicity $\lambda\mu$ in $G(\lambda)\times H(\mu)$ and $G(\lambda)\times H(\mu)$ is isomorphic to $(G\times H)(\lambda\mu).$ Hence $G(\lambda)\times H\cong G\times H(\lambda)\cong (G\times H)(\lambda).$

If $H_1,H_2,\ldots,H_k$ are edge-disjoint subgraphs of $G$ and $E(G)=E(H_1)\cup E(H_2)\cup\ldots\cup E(H_k),$ then we write $G=H_1\oplus H_2\oplus\ldots\oplus H_k.$ For a graph $G,$ if $E(G)$ can be partitioned into $E_1,$$E_2,\ldots,$$E_k$ such that $\langle E_i\rangle\cong H,$ for all $i,\,1\le i\le k,$ then we say that $H$ decomposes $G,$ or {\it $H$-decomposition} of $G$ exists. Clearly, the tensor product is commutative and distributive over edge-disjoint union of graphs, that is, if $G=H_1\oplus H_2\oplus \ldots \oplus H_k,$ then $G\times H=(H_1\times H)\oplus (H_2\times H)\oplus\ldots\oplus(H_k\times H).$

Let $G$ be a finite group and let $S$ be a symmetric subset of $G$ (that is, $s\in S \text{ implies } -s\in S).$ The vertices of the {\it Cayley graph,} $Cay(S,\, G),$ are the elements of $G$ and there is an edge between $x$ and $y$ if and only if $x-y \in S.$ Note that $Cay(S,\, G)$ is connected if and only if $S$ generates the group $G.$ A \emph{circulant} $X=Circ(n;\,L)$ is a graph with vertex set $V(X)=\{u_0,\,u_1,\,\ldots,\,u_{n-1}\}$ and edge set $E(X)=\{u_iu_{i+\ell}\,|\,i\in \mathbb{Z}_n,\,\ell \in L\},$ where $L\subseteq\left\{1,\,2,\,\ldots,\,\lfloor\frac{n}{2}\rfloor\right\}$ and $\mathbb{Z}_n$ is the set of integers modulo $n.$ The elements of $L$ are called {\it jumps.} Clearly, every circulant graph of order $n$ is a Cayley graph with the underlying group being $\mathbb{Z}_{n}.$

The problem of finding Hamilton cycle decompositions of product graphs is not new. Hamilton cycle decompositions of various product graphs, including digraphs, have been studied by many authors; see, for example, \cite {{CyclesandRays}, {Ann.DM.3.21-28}, {Bosak.Book}, {IJPAM.23.723-729}, {Comput.Mat.Appl.31.11-19}, {ARS.80.33-44},  {DM.308.3586-3606}, {DM.310.2776-2789}, {JCTB.73.119-129}, {DM.90.169-190}, GC.25.571-581}. It has been conjectured \cite{Ann.DM.3.21-28} that if both $G$ and $H$ are Hamilton cycle decomposable graphs, then $G\Box H$ is Hamilton cycle decomposable, where $\Box$ denotes the cartesian product of graphs. This conjecture has been verified to be true for a large classes of graphs \cite{DM.90.169-190}. Baranyai and Sz\'{a}sz \cite{JCTB.31.253-261} proved that if both $G$ and $H$ are even regular Hamilton cycle decomposable graphs, then $G \circ H$ is Hamilton cycle decomposable. In \cite{JCTB.73.119-129}, Ng has obtained a partial solution to the following conjecture of Alspach et al. \cite{CyclesandRays}: If $D_1$ and $D_2$ are directed Hamilton cycle decomposable digraphs, then $D_1\circ D_2$ is directed Hamilton cycle decomposable. Jha \cite{IJPAM.23.723-729} conjectured the following

\begin{con}\emph{\cite{IJPAM.23.723-729}}\label{IJPAM.23.723-729}
If both $G$ and $H$ are Hamilton cycle decomposable graphs and $G\times H$ is connected, then $G\times H$ is Hamilton cycle decomposable.
\end{con}

Conjecture \ref{IJPAM.23.723-729} was disproved, see \cite{DM.186.1-13}. Because of this, finding Hamilton cycle decompositions of the tensor products of Hamilton cycle decomposable graphs is considered to be difficult. Though Conjecture \ref{IJPAM.23.723-729} has been disproved, we believe that if the graphs $G$ and $H$ are suitably chosen, that is, with some suitable conditions imposed on them, then $G\times H$ may have Hamilton cycle decomposition. In \cite{DM.268.49-58} and \cite{ARS.80.33-44} it has been proved that $K_r \times K_s$ and $K_{r ,r} \times K_m$ are Hamilton cycle decomposable; in \cite{DM.308.3586-3606} it is shown that the tensor product of two regular complete multipartite graphs is Hamilton cycle decomposable. Hamilton cycle decompositions of the tensor products of complete bipartite graphs and complete multipartite graphs are dealt with in \cite{DM.310.2776-2789}. Also in \cite{GC.25.571-581}, Paulraja and Sivasankar proved that $(K_r\times K_s)^*,\,((K_r \circ \overline{K}_s)\times K_n)^*,\, ((K_r \times K_s)\times K_m)^*,\,((K_r\circ \overline{K}_s)\times (K_m\circ \overline{K}_n))^*$ and $(K_{r,r} \times (K_m\circ \overline{K}_n))^*$ are directed Hamilton cycle decomposable. It can be observed that $K_r,\,K_{r,\,r},\,K_r\circ\overline{K}_s$ are circulant graphs. Based on the results of \cite{{IJPAM.23.723-729}, {Comput.Mat.Appl.31.11-19}, {ARS.80.33-44},  {DM.308.3586-3606}, {DM.310.2776-2789}}, Manikandan and Paulraja conjectured the following:

\begin{con}\emph{(Manikandan and Paulraja) \cite{Thesis.RSM}.}\label{con.RSMPP}
If $G$ and $H$ are Hamilton cycle decomposable circulant graphs and at least one of them is non bipartite, then $G\times H$ is Hamilton cycle decomposable.
\end{con}

\begin{thm}\emph{\cite{Ann.DM.3.21-28}}\label{Ann.DM.3.21-28}
If both $G$ and $H$ have Hamilton cycle decompositions and at least one of them is of odd order, then $G\times H$ admits a Hamilton cycle decomposition.\hfill$\Box$
\end{thm}

One may naturally ask when is $G\times H$ Hamilton cycle decomposable, if both $G$ and $H$ are of even order. In this paper it is partially answered.

 We say that an even regular circulant graph $X=Circ(n;\,L)$ has the {\it property $Q,$} if (1) the number of odd and even jumps in $L$ are equal; (2) odd jumps can be paired with even jumps so that each of the $\frac{|L|}{2}$ resulting $4$-regular graphs is connected. It is known that every $4$-regular connected circulant graph is Hamilton cycle decomposable; see \cite{JCTB.46.142-153}. Consequently, every circulant graph with property $Q$ is Hamilton cycle decomposable.

Here we prove the following main Theorem.
\begin{thm}\label{THM.CIR.1}
Let $G$ be a circulant graph with property $Q$ and let $H$ be any Hamilton cycle decomposable multigraph, then $G\times H$ is Hamilton cycle decomposable.
\end{thm}


This theorem has many interesting consequences. In particular, we have the following corollary.
\begin{cor}\label{COR.CIR.1}
If $G$ and $H$ are even regular Hamilton cycle decomposable circulant graphs and at least one of them has the property $Q,$ then $G\times H$ is Hamilton cycle decomposable.
\end{cor}

One of the consequences of Corollary \ref{COR.CIR.1} is that if $G=(K_{4n+2}-F)$ and $H=(K_{2m}-F^\prime),$ where $F$ and $F^\prime$ are $1$-factors of $K_{4n+2}$ and $K_{2m},$ respectively, then $G\times H$ is Hamilton cycle decomposable. In particular, after deleting suitable number of jumps of $K_{4n+2}$ and $K_{2m},$ the resulting graphs need not be dense but their tensor product is Hamilton cycle decomposable. This cannot be deduced from the existing results in this direction.


\section{\bf Notation and Preliminaries.}

First we present the necessary definitions here. The notation that we use here are from \cite{JCTB.46.142-153} and for the sake of completeness we give them.

The Cayley graph $\Gamma=Cay(n,\,S),$ where $S=\{a,\,b\}$ is a generating set of the finite abelian group with $2a\ne 0,\,2b\ne 0$ and $a\ne\pm b,$ is a simple connected graph. We call the edge $x(x+a)$ of $\Gamma$ an $a$-edge; the subgraph formed by the $a$-edges is a disjoint union of cycles, called $a$-cycles, each of length $k_a,$ the order of the element $a.$ Similarly, we define $b$-edges and $b$-cycles; each $b$-cycle is of length $k_b,$ the order of the element $b.$ In $\Gamma,$ let us denote the number of $a$-cycles by $\alpha$ and the number of $b$-cycles by $\beta.$ Since the length of each $a$-cycle (respectively $b$-cycle) is $k_a=\frac{n}{\alpha},$ (respectively $k_b=\frac{n}{\beta},$) $n=\alpha k_a$ (respectively $n=\beta k_b$).

As the graph $\Gamma$ is connected, the vertices $0,\,b,\,\ldots,\,\ell b,\,\ldots,\,(\alpha-1)b$ are in the $\alpha$ $a$-cycles denoted by $C_0,\,C_1,\,\ldots,C_\ell,\,\ldots,C_{\alpha-1},$ respectively, and $\alpha b$ belongs to $C_0,$ see \cite{JCTB.46.142-153}.  Hence, we have $\alpha b=ca$ for some $c$ with $0\le c\le k_a-1.$ Further, every vertex $x$ of $\Gamma$ can be uniquely written as $x=ib+ja$ with $i\in \{0,\,1,\,\ldots,\, \alpha-1\}$ and $j\in \{0,\,1,\,\ldots,\, k_a-1\}.$ Using this uniqueness, label the vertices of $\Gamma$ as $(i,\,j)$ with $i\in \{0,\,1,\,\ldots,\, \alpha-1\},$ and $j\in \{0,\,1,\,\ldots,\, k_a-1\},$ where the first coordinate indicates the label of the cycle $C_i$ containing the vertex and the second coordinate indicates the position of the vertex on the cycle.

The following definition comprises the notation and defines a class of simple graphs, which we denote by $\Gamma(\alpha,\,\beta).$

\begin{defn}\emph{\cite{JCTB.46.142-153}}\label{JCTB.46.142-153}
Let $\Gamma(\alpha,\,\beta)$ denote the class of simple graphs on $\alpha k$ vertices, where $\alpha\ge 1,\,k\ge 3,\,0\le c<k,$ and $\beta=gcd(k,\,c).$ The $\alpha k$ vertices of the graph can be labeled $(i,\,j)$ with $i$ taken modulo $\alpha$ and $j$ taken modulo $k.$ The edges are \emph{(1)} First kind: $(i,\,j)(i,\,j+1)$ and \emph{(2)} Second kind: $(i,\,j)(i+1,\,j)$ for all $i\in \{0,\,1,\,\ldots,\, \alpha-2\},$ and $(\alpha-1,\,j)(0,\,j+c).$
\end{defn}

A graph $\Gamma\in\Gamma(\alpha,\,\beta)$ has $\alpha$ \lq\lq vertical" disjoint cycles $C_i,\,0\le i<\alpha,$ with a natural orientation and $\alpha-1$ horizontal {\it parallel matchings} between the cycles $C_i$ and $C_{i+1}$ for $0\le i<\alpha-1,$ and a particular {\it parallel matching} between $C_{\alpha-1}$ and $C_0$ (which depends on the value of $c).$ A graph $\Gamma\in\Gamma(2,\,\beta)$ consists of two cycles plus two perfect matchings between them. A graph $\Gamma\in\Gamma(1,\,\beta)$ consists of a cycle plus the chords joining $(0,\,j)$ to $(0,\,j+c).$ Observe that $\Gamma(\alpha,\,\beta)$ and $\Gamma(\beta,\,\alpha)$ are isomorphic classes of graphs, that is, an element $\Gamma_1\in\Gamma(\alpha,\,\beta)$ is isomorphic to an element $\Gamma_2\in\Gamma(\beta,\,\alpha)$ and vice versa, for example see Figure 1.

\begin{center}
\includegraphics{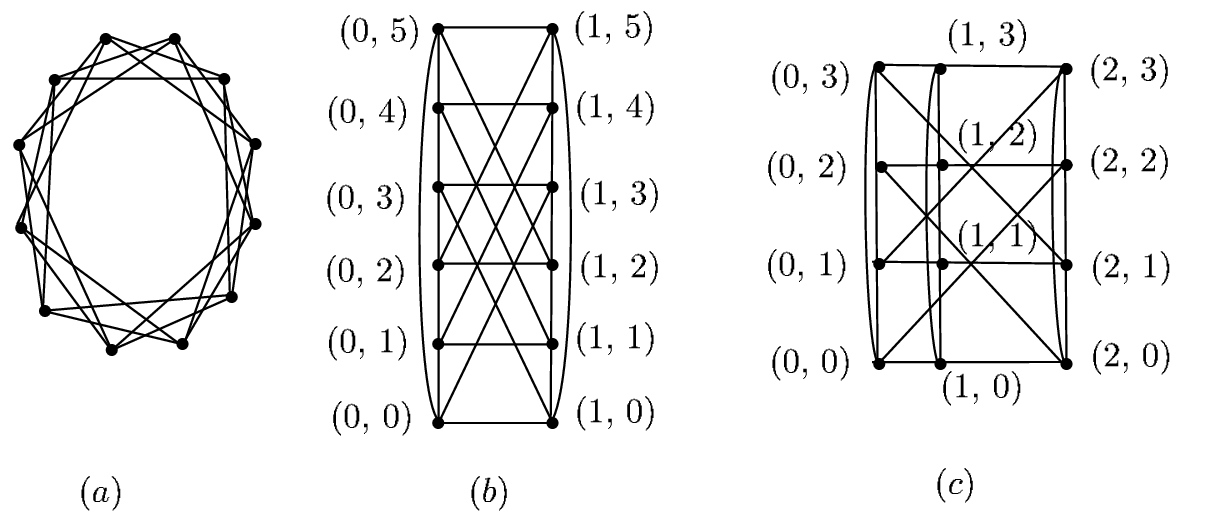}
\end{center}
{\small Graphs of Figures $1(a),\,1(b)$ and $1(c)$ are isomorphic. \\$(a)$ Graph $G=Circ(12,\,\{2,\,3\}),$ where the edges of jump $2$ give $2$ disjoint cycles each of length $6$ and the edges of jump $3$ give $3$ disjoint cycles each of length $4.$  \\$(b)$ $G$ is an element of $\Gamma(2,\,3)$ with $c=3.$ \\$(c)$ $G$ in $(b)$ is drawn as an element of $\Gamma(3,\,2)$ with $c=2.$
}\begin{center}Figure 1\end{center}

\begin{defn}\emph{\cite{JCTB.46.142-153}}\label{JCTB.46.142-153-defn}
A Hamilton cycle decomposition of a graph $\Gamma\in\Gamma(\alpha,\,\beta)$ has the property $Q_1,$ if both Hamilton cycles use at least one edge of the matching between $C_i$ and $C_{i+1},$ for all $i\in\{0,\,1,\,\ldots,\,\alpha-1\},$ where $C_\alpha=C_0.$
\end{defn}

\begin{thm}\emph{\cite{JCTB.46.142-153}}\label{JCTB.46.142-153-1}
The class $\Gamma(\alpha,\,\beta)$ consists of the $4$-regular connected Cayley graphs on a finite abelian group with a generating set $\{a,\,b\},$ where $\alpha$ is the number of $a$-cycles and $\beta$ is the number of $b$-cycles. \hfill$\Box$
\end{thm}

Using the above theorem, Bermond et al. proved the following
\begin{thm}\emph{\cite{JCTB.46.142-153}}\label{JCTB.46.142-153-main}
Every $4$-regular connected Cayley graph on a finite abelian group can be decomposed into two Hamilton cycles.
\end{thm}

If $G$ and $H$ are Hamilton cycle decomposable graphs with at least one of them having odd order, then $G\times H$ is Hamilton cycle decomposable, by Theorem \ref{Ann.DM.3.21-28}; hence in what follows, {\bf we assume that both $G$ and $H$ are Hamilton cycle decomposable graphs each having an even number of vertices.} If $G$ is a connected $4$-regular circulant graph of even order with generating set $\{a,\,b\}$ where $a$ and $b$ are of different parity such that $2a\ne 0,\,2b\ne 0$ and $a\ne\pm b,$ then by the Definition \ref{JCTB.46.142-153}, $\alpha$ and $\beta$ are of different parity. For our convenience, according to the situation we consider the parity of $a$ and $b.$ {\bf In what follows, we assume that $\Gamma(\alpha,\,\beta)$ denotes the class of Cayley graphs on the abelian group $\mathbb{Z}_n$} (that is, circulants) with generating set $\{a,\,b\}$ and by Theorem \ref{JCTB.46.142-153-1}, $\Gamma(\alpha,\,\beta)$ consists of the class of $4$-regular connected circulants of order $n$ with generating set $\{a,\,b\},$ where $\alpha$ is the number of $a$-cycles in a $2$-factor and $\beta$ is the number of $b$-cycles.

In the rest of the paper, {\bf we assume that $\Gamma\in\Gamma(\alpha,\,\beta)$ is a $4$-regular circulant graph of even order with generating set $\{a,\,b\}$ where $a$ and $b$ are of different parity; consequently, $\alpha$ and $\beta$ are of different parity.}

\begin{defn}\label{defn.jha}
Let $\{H_1,\,H_2\}$ be a Hamilton cycle decomposition of a $4$-regular graph $G$ of even order with at least $6$ vertices.  If $G$ contains a $4$-cycle $(a\,b\,c\,d)$ such that the edges $ab$ and $cd$ belong to one of the two Hamilton cycles and the edges $bc$ and $da$ are on the other Hamilton cycle, then the $4$-cycle is said to be an {\it alternating} $4$-cycle in $G,$ with respect to $H_1$ and $H_2.$ Further, if the vertices $a$ and $c$ \emph{(}or $b$ and $d$\emph{)} are at an odd distance along each of the two Hamilton cycles $H_i,$ then the $4$-cycle $(a\,b\,c\,d)$ is said to be an odd alternating $4$-cycle in $G$ with respect to $\{H_1,\,H_2\}.$ We say that the graph $G$ has {\it property $Q_2,$} if $G$ contains an odd alternating $4$-cycle with respect to a Hamilton cycle decomposition of it.
\end{defn}

The proof techniques we use here heavily depend on \cite{JCTB.46.142-153}.

The {\it reduced graph} $\Gamma^\prime\in\Gamma(\alpha,\,\beta)$ with respect to $\Gamma\in\Gamma(\alpha+2,\,\beta)$ is defined as follows: delete the vertices of the cycles $C_\alpha$ and $C_{\alpha+1}$ of $\Gamma,$ that is, the vertices $(\alpha,\,j)$ and $(\alpha+1,\,j),\,0\le j\le k-1,$ and add the edges $(\alpha-1,\,j)(0,\,j+c),\,0\le j\le k-1,$ joining $C_{\alpha-1}$ and $C_0,$ where $c$ is given in the definition of $\Gamma(\alpha+2,\,\beta).$ In some cases $\Gamma^\prime$ might be a multigraph and, if it is a simple graph, it is an element of $\Gamma(\alpha,\,\beta)$ ($k$ and $c$ being unchanged). Successive reductions of the resulting graphs result in a graph $\Gamma^{k}\in\Gamma(\alpha+2-2k,\,\beta)$ for some $k\ge 1.$ $\Gamma^{k}\in\Gamma(\alpha+2-2k,\,\beta)$ is also called a reduced graph of $\Gamma\in\Gamma(\alpha+2,\,\beta).$

The graph $\Gamma\in\Gamma(\nu+2,\,\delta)$ is said to be {\it lift graph} of  $\Gamma^\prime\in\Gamma(\nu,\,\delta),$ if $\Gamma$ is obtained as follows: delete the edges  $(\nu-1,\,j)(0,\,j+c),\,0\le j\le k-1,$ joining $C_{\nu-1}$ and $C_0$ of $\Gamma^\prime;$ then add two $a$-cycles $C_{\nu}$ and $C_{\nu+1}$ with vertices $(\nu,\,j)$ and $(\nu+1,\,j),\,0\le j\le k-1,$ to $\Gamma^\prime,$ and add the edges  $(\nu-1,\,j)(\nu,\,j),\,(\nu,\,j)(\nu+1,\,j),\,(\nu+1,\,j)(0,\,j+c),\,0\le j\le k-1,$ where $c$ is as given in the definition of $\Gamma(\nu,\,\delta).$ Successive liftings of the resulting graphs result in a graph  $\Gamma^{k}\in\Gamma(\nu+2k,\,\delta)$ for some $k\ge 1.$ $\Gamma^{k}\in\Gamma(\nu+2k,\,\delta)$ is also called a lift graph of $\Gamma^\prime\in\Gamma(\nu,\,\delta).$

We use $p_G(x,\,y)$ to denote the length of a path (not necessarily a shortest path) from $x$ to $y$ in the graph $G.$ By $(a,\,b)$-section, we denote a path from $a$ to $b.$

\begin{lem}\label{LEM2.directcirculant}
Every graph $\Gamma\in\Gamma(\alpha,\,\beta),$ where $1\le \alpha,\,\beta\le 2,$ has a Hamilton cycle decomposition $\{H_1,\,H_2\}$ with the properties $Q_1$ and $Q_2.$
\end{lem}

{\bf Proof.} Let $\Gamma\in\Gamma(\alpha,\,\beta).$ Since $1\le \alpha,\,\beta\le 2$ and $\alpha\ne \beta,$ without loss of generality assume that $\alpha=2$ and $\beta=1.$ Consequently, $a$ is even and $b$ is odd as $2=\alpha=gcd(n,\,a)$ and $1=\beta=gcd(n,\,b).$ Clearly, $c\ne 0$ (otherwise $\beta\ne 1).$ The existence of Hamilton cycles $H_1$ and $H_2$ in $\Gamma$ described below is in \cite{JCTB.46.142-153}.
A Hamilton cycle $H_{1}$ is obtained by deleting the edges $0a$ and $(-b)(-b+a)$ from $C_0$ and $C_1,$ respectively, and adding the edges $(-b)0$ and $(-b+a)a$ connecting these two $a$-cycles, that is, $H_1=\{C_0-\{0a\}\}\cup \{C_1-\{(-b)(-b+a)\}\}\cup \{0(-b),\,(-b+a)a\},$  see Figure 2(b);  $H_2=\Gamma-E(H_1)$ is shown in the solid lines of Figure 2(c1). Clearly, these two Hamilton cycles use both $a$ and $b$-edges and hence the property $Q_1$ holds.

\begin{center}
\includegraphics[scale=0.8]{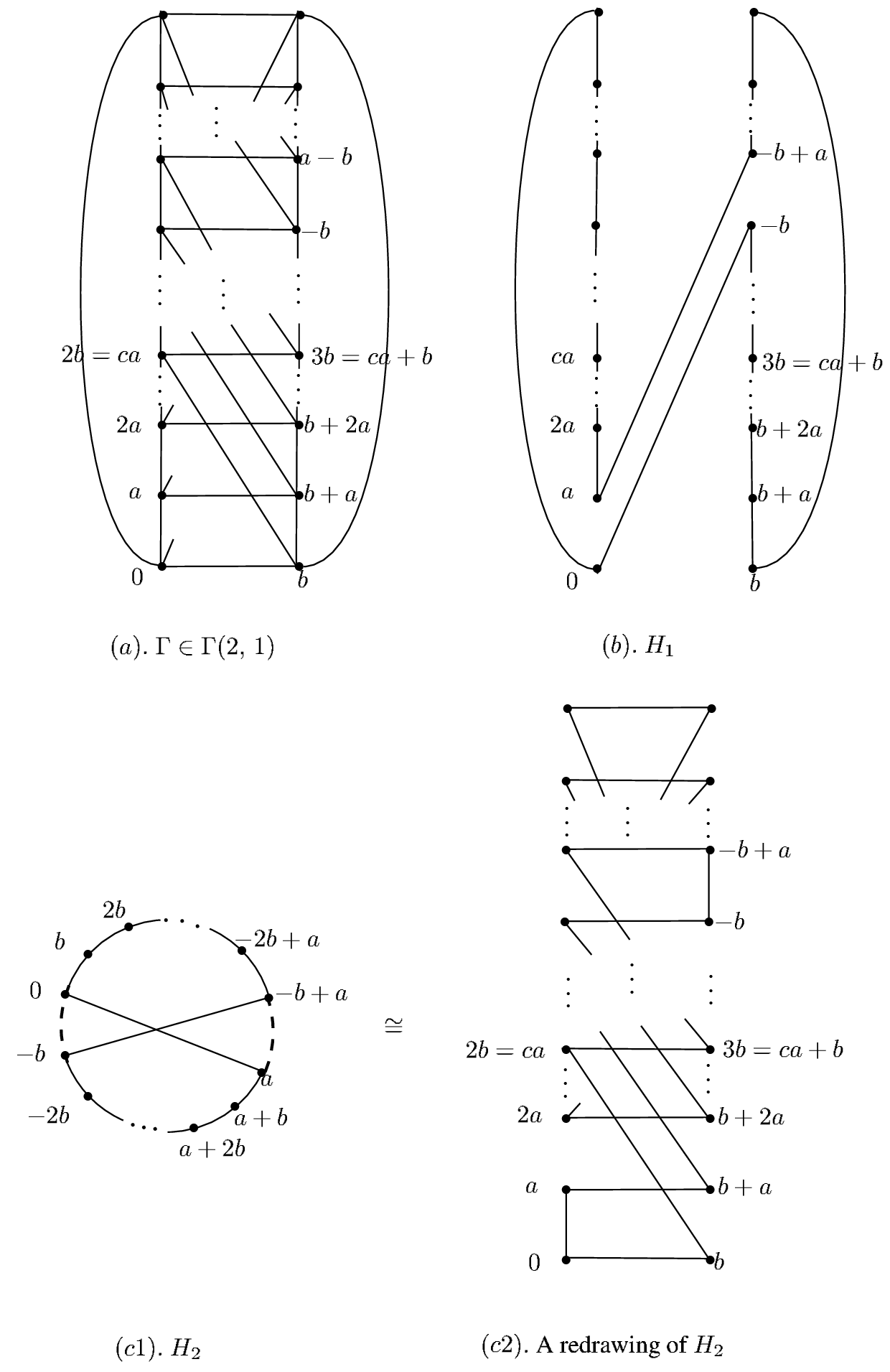}
\end{center}
{\small (a). A member $\Gamma$ of $\Gamma(2,\,1)$\\ (b). Hamilton cycle $H_1$ of $\Gamma\in\Gamma(2,\,1)$\\ (c1). The solid lines are edges on $H_2.$ Other than the two $a$-edges $0a$ and $(-b)(-b+a),$ all edges of $H_2$ (and thus all edges connecting $C_0$ and $C_1)$ are $b$-edges. \\ (c2). Another drawing of the Hamilton cycle $H_2$ of $\Gamma\in\Gamma(2,\,1).$}\begin{center}Figure 2
\end{center}

\NI{\bf Claim.} The graph $\Gamma\in\Gamma(2,\,1)$ satisfies the property $Q_2$ with respect to the above Hamilton cycle decomposition $\{H_1,\,H_2\}.$

Since $\alpha=2,$ there are two $a$-cycles $C_0$ and $C_1$ and there are two parallel matchings between $C_0$ and $C_1.$ We denote one of the matchings from $C_0$ to $C_1$ as $M_0$ and the other matching from $C_1$ to $C_0$ as $M_1;$ let the edges of $M_0$ be $(0,\,j)(1,\,j),\,0\le j\le \frac{n}{2}-1,$ and let the edges of $M_1$ be $(1,\,j)(0,\,j+c),\,0\le j\le \frac{n}{2}-1,$ where addition is taken modulo $n.$ We call the vertex $(i,\,j)\,(=ib+ja)$ of $C_0$ as the {\it corresponding vertex} of $(i+1,\,j)\,(=(i+1)b+ja)$ of $C_1$ and vice versa. We prove the existence of the property $Q_2$ with respect to the Hamilton cycle decomposition $\{H_1,\,H_2\}$ of $\Gamma$ by the parity of $c.$

\NI{\bf Case 1.} $n\equiv 2\,(mod\,\,4).$

First, we shall obtain an even length path along $H_1,$ so that this path is a section of a suitable odd length path along $H_1$ having its ends as two \lq\lq opposite vertices" of an odd alternating $4$-cycle in $\Gamma;$ the existence of the odd alternating $4$-cycle will be proved later.

In ${\bf (A)}$ below, we obtain a required even length path along $H_1$ and in ${\bf (B)}$ we prove that the path obtained in ${\bf (A)}$ is a section of a path along $H_1$ joining two \lq\lq opposite vertices" of an odd alternating $4$-cycle.

\NI{\bf (A).} First we consider $c$ to be even and consider the Hamilton cycle $H_1$ of $\Gamma$ obtained above. Since $M_1$ is a perfect matching of $\Gamma$ and $M_1$ matches the vertices in $C_1$ with vertices in $C_0,$ $(-b)0$ is an edge of $M_1,$ which is one of the two $M_1$ edges in $H_1,$ see Figure 3(a). Clearly, $-b\equiv (\frac{n}{2}-c)a+b\,(mod\,\,n),$ that is, $-b$ is the $(\frac{n}{2}-c+1)$th vertex of $C_1,$ starting from $b,$ see Figure 3(a). As $\frac{n}{2}$ is odd, $(\frac{n}{2}-c+1)$ is even and hence the length of the section of $C_1$ from $b$ to $-b$ (containing the vertex $b+a$) is odd.
\begin{center}
\includegraphics[scale=0.8]{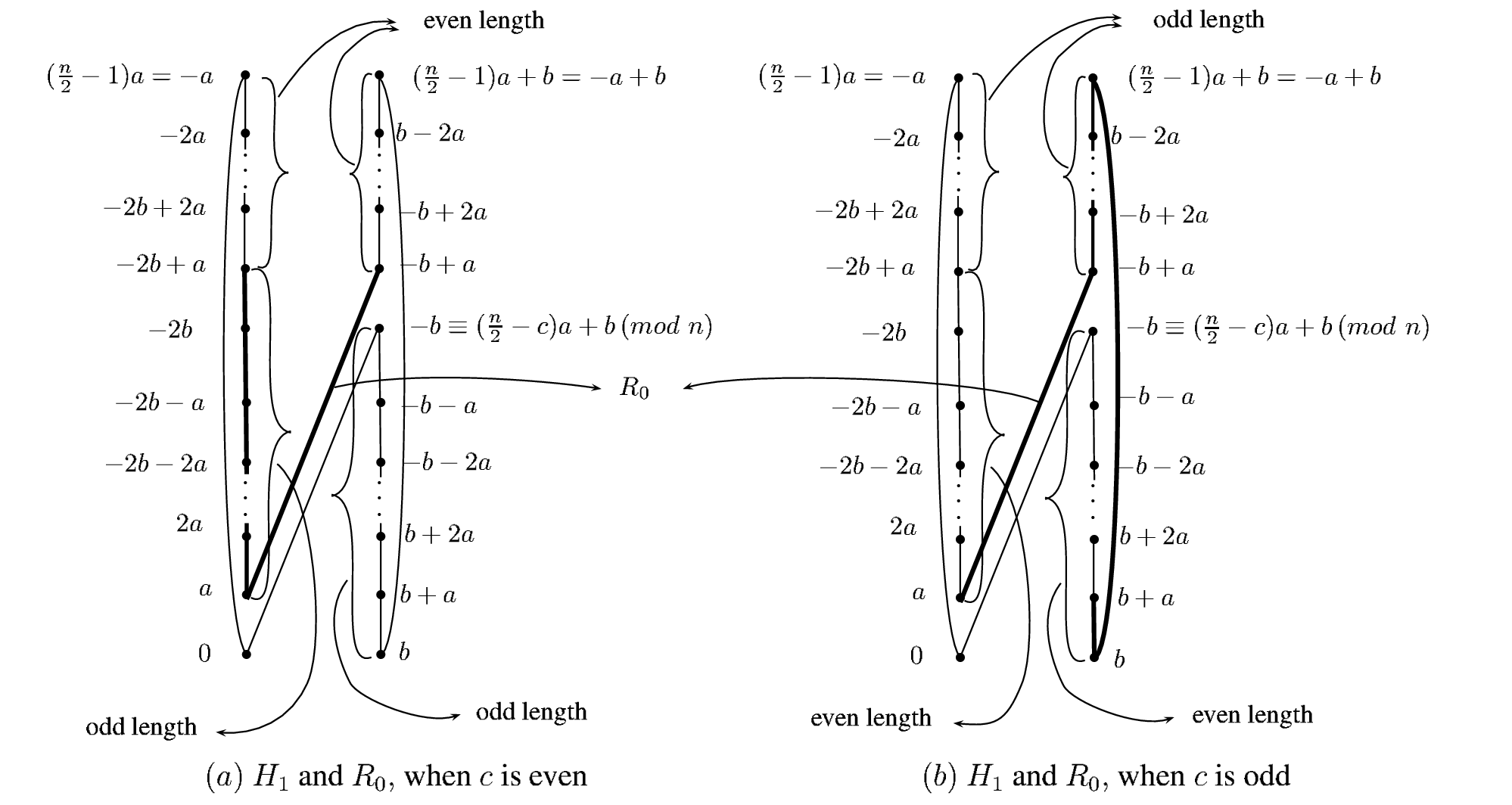}
{\small $H_1$ is shown in both Figures 1(a) and 1(b), wherein the bold lines are the required even length path.}\\Figure 3
\end{center}

As the length of the section of $C_1$ from $b$ to $-b$ (containing the vertex $(b+a)$) is the same as the length of the section of $C_0$ from $a$ to $(-2b+a)$ (containing the vertex $2a),$ the length of the path from $a$ to $-2b+a$ along $H_1,$ denoted by $p_{H_1}(a,\,-2b+a)$ is odd and let this path be $R.$ Now $R$ together with the edge $a(-b+a)$ is our required even length path $R_0,$ see the bold edges of Figure 3(a).

Next we suppose that $c$ is odd. As in the above paragraph, the vertex $0$ in $C_0$ is matched under $M_1$ with the vertex $-b\equiv (\frac{n}{2}-c)a+b\,(mod\,\,n);$ $\frac{n}{2}-c$ is even as $\frac{n}{2}$ is odd. Thus the length of the section of $C_1$ from $b$ to $-b$ (containing the vertex $(b+a)$) is even and hence the length of the section $R_1$ of $C_1$ from $b$ to $(-b+a)$ (containing the vertex $(-a+b)$) is even (note that $(-b)(-b+a)$ is an edge of $C_1$ which we have deleted for the construction of $H_1$). Now $R_1$ together with the edges $\{(b+a)b,\,(-b+a)a\}$ is our required even length path $R_0$ from the vertex $a$ to the vertex $b+a$ along $H_1,$ see the bold edges of Figure 3(b).

\NI ${\bf (B).}$ $H_1$ contains exactly two edges of $M_1$ as shown in the Figure 3. Hence except the two $4$-cycles of $\Gamma,$ namely, $(0,\,a,\,a+b,\,b)$ and $(-2b,\,-2b+a,\,-b+a,\,-b),$ the other $4$-cycles, which have two consecutive vertices of $C_0$ and the corresponding two vertices of $C_1$ constitute alternating $4$-cycles (with respect to $H_1$ and $H_2$ of $\Gamma$).

If $c$ is even, it is clear that every $4$-cycle containing two consecutive vertices of the section $(-2b+a)(-2b+2a)(-2b+3a)\ldots((\frac{n}{2}-1)a)0$ of $C_0$ and the two corresponding consecutive vertices of the section $(-b+a)(-b+2a)(-b+3a)\ldots(b-a)b$ of $C_1$ is an alternating $4$-cycle of $\Gamma$ with respect to $H_1$ and $H_2,$ see the Figure 3 (note that $a=2b$ is ruled out, for otherwise, it would imply $c=1$ for the following reason: $a=2b$ implies $b=-b+a$ and since $(-b+a)a$ is an edge of $H_1$ joining $C_1$ and $C_0$ implies $c=1,$ which is not the case as we consider $c$ is even). Similarly if $c$ is odd, every $4$-cycle containing two consecutive vertices in the section $a\,(2a)\,(3a)\,\ldots\,(-2b)$ of $C_0$ and the two corresponding consecutive vertices in the section $(b+a)\,(b+2a)\,(b+3a)\,\ldots\,(-b)$ of $C_1$ is an alternating $4$-cycle of $\Gamma$ with respect to $H_1$ and $H_2;$ note if $a=-2b,$ then $a=-ca$ so $-a=ca=(\frac{n}{2}-1)a$ but then $c\equiv (\frac{n}{2}-1)\,(mod\,\,\frac{n}{2}),$ that is, $c=\frac{n}{2}-1,$ which is even, a contradiction.

Next we prove that the alternating $4$-cycles described in the above paragraph using the sections of $C_0$ and $C_1$ satisfy the property that the length of the path along $H_1$ between the opposite vertices of the $4$-cycle is of odd length and further this odd length path contains $R_0$ (described above).

Let $C=(x,\,x+a,\,x+a+b,\,x+b)$ be an alternating $4$-cycle so that $x,\,x+a\in V(C_0)$ and $x+a+b,\,x+b\in V(C_1),$ see Figure 4. The vertices $x$ and $x+b$ are the corresponding vertices in $C_0$ and $C_1.$ If $c$ is even, then the $(-2b+a,\,x)$-section of $H_1,$ contained in $C_0,$ and the corresponding $(x+b,\,(-b+a))$-section of $H_1$ contained in $C_1$ have the same length, see Figure 4(a). Hence the $(x,\,x+a+b)$-section of $H_1$ (containing the vertex $-2b$) is of odd length as it contains $R_0,$ which is of even length, see Figure 4(a); that is $p_{H_1}(x,\,x+a+b)$ is odd when $c$ is even. Similarly, if $c$ is odd then, the $(a,\,x)$-section of $H_1$ (containing the vertex $2a$), contained in $C_0,$ and the corresponding $(b+a,\,x+b)$-section of $H_1,$ contained in $C_1$ (containing the vertex $b+2a),$ have the same length. Thus the $(x,\,x+a+b)$-section of $H_1$ (containing the vertex $-b+a)$ is of odd length as it contains $R_0,$ which is of even length, see Figure 4(b). That is $p_{H_1}(x,\,x+a+b)$ is odd when $c$ is odd, see Figure 4(b).

\begin{center}
\includegraphics[scale=0.8]{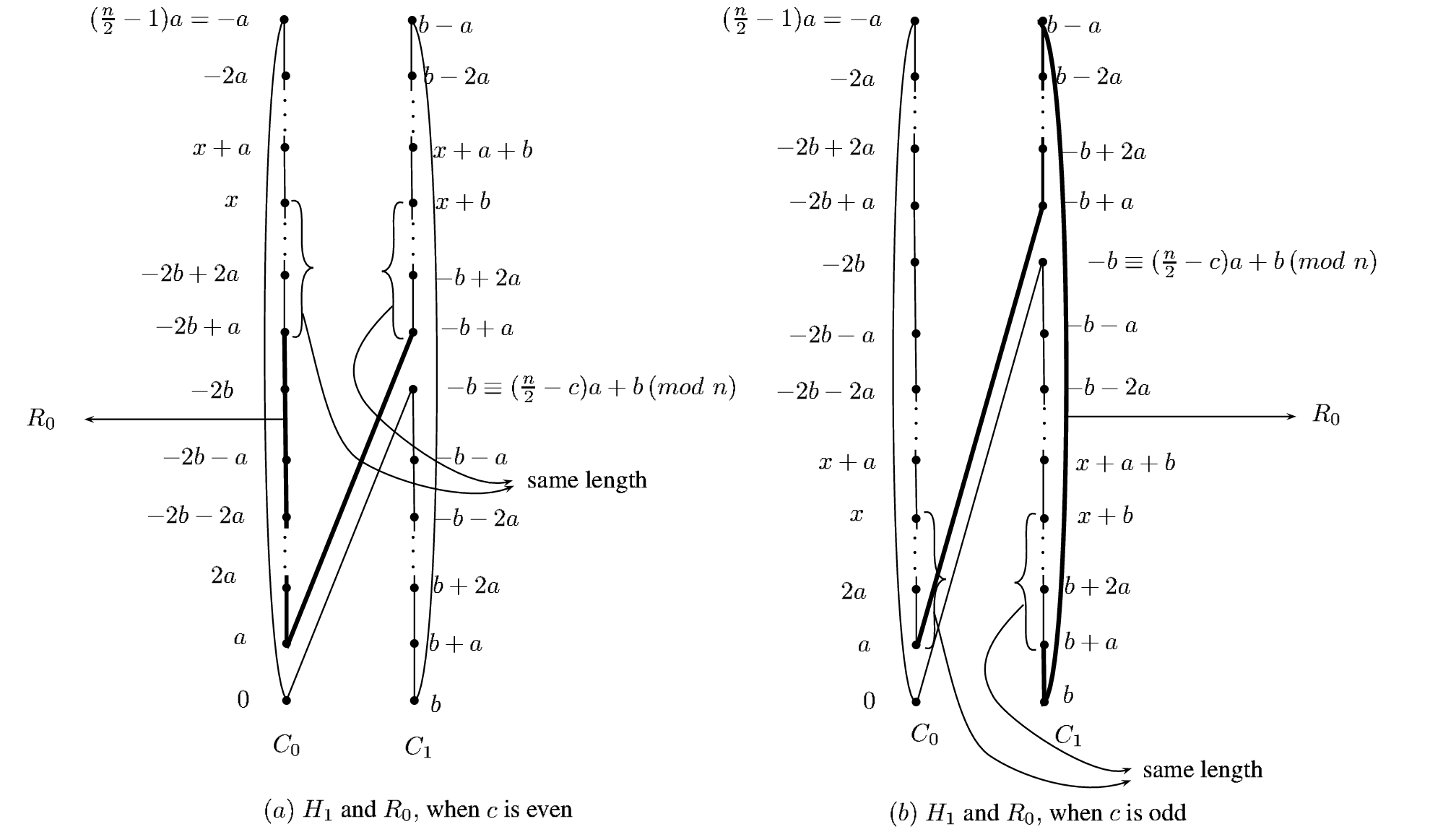}
\end{center}
{\small (a) The lengths of the sections $(-2b+a,\,x)$ along $H_1$ (section of $C_0$) and $(-b+a,\,x+b)$ along $H_1$ (section of $C_1$) are same and $(x+b)(x+a+b)$ is an edge \\(b) the lengths of the sections $(a,\,x)$ along $H_1$ (part of $C_0$) and $(b+a,\,x+b)$ along $H_1$ (part of $C_1$) are same and $(x+a)(x+a+b)$ is an edge.}\begin{center}Figure 4
\end{center}

Next we shall obtain an appropriate $x$ on $C_0$ so that $p_{H_2}(x,\,x+a+b)$ is also odd which gives the property $Q_2$ with respect to $H_1$ and $H_2.$ As $H_2$ contains exactly two $a$-edges, namely, $(-b)(-b+a)$ and $0a$ (see Figure 2(c2)), $H_2-\{0a,\,(-b)(-b+a)\}$ is a pair of odd length paths (since the vertices of these paths alternate between $V(C_0)$ and $V(C_1)$). Let $S_0=0\,(b)\,(2b)\,\ldots\,(-b+a)$ and $S_1=a\,(a+b)\,(a+2b)\,\ldots\,(-b)$ denote these paths, namely, the $(0,\,(-b+a))$-section and $(-b,\,a)$-section of $H_2$, respectively, see Figures 2(c1) and 2(c2).

If $S_0$ (respectively $S_1$) contains a pair of consecutive vertices $x$ and $x+a$ of $C_0$ (note that the edge $x(x+a)$ is in $H_1$ and $S_0$ is contained in $H_2$), then $p_{H_2}(x,\,x+a+b)$ is odd as the vertices of $S_0$ (respectively $S_1$) alternate between $V(C_0)$ and $V(C_1)$ and $(x+a)(x+a+b)$ is an edge in $H_2.$ We shall prove the existence of such a pair of consecutive vertices in the sequence of vertices $L_0=(-2b+a)\,(-2b+2a)\,(-2b+3a)\,\ldots\,(\frac{n}{2}-1)a\,0,$ or in $L_1=a\,(2a)\,(3a)\,\ldots\,(-2b),$ see Figure 4.

Let $c$ be even. From the construction of $H_1$ and $H_2,$ the vertices $0$ and $-2b+a$ are in $S_0,$ (see Figure 2(c1)). There are even number of vertices in $L_0$ as $R_0$ is of even length and $C_0$ has odd number of vertices (see Figure 3(a)). If there is no pair of consecutive vertices in $L_0$ of the required type in $S_0$ or $S_1,$ then the preceding vertex of $0,$ namely, $(\frac{n}{2}-1)a\,(=-a),$ and the succeeding vertex of $-2b+a,$ namely, $-2b+2a,$ are in $S_1,$ as the vertices $0$ and $-2b+a$ are in $S_0;$ consequently, there must be an odd number of vertices from $-2b+2a$ to $(\frac{n}{2}-1)a(=-a)$ in $L_0,$ which is not the case (since in $L_0,$ $-2b+2a$ to $-a$ contains even number of vertices). Therefore there must exist a pair of consecutive vertices in $L_0$ of the required type and hence $\Gamma$ satisfies property $Q_2$ with respect to $H_1$ and $H_2.$

Next we assume that $c$ is odd. As above, we shall show the existence of the required pair of consecutive vertices in $L_1=a\,(2a)\,(3a)\,\ldots\,(-2b).$ $L_1$ has even number of vertices (see Figure 3(b)). Assume that there is no pair of vertices in $L_1$ of the required type in $S_0$ or $S_1$ and hence alternate vertices of $L_1$ are in $S_0$ and $S_1.$ Then, as the vertices $a$ and $-2b$ are in $S_1,$ the vertices $2a$ and $-2b-a$ must be in $S_0$ and hence there must be an odd number of vertices from $2a$ to $-2b-a$ in $L_1,$ which is not the case (as $L_1$ contains even number of vertices from $2a$ to $-2b-a$). Thus there must exist a pair of consecutive vertices as required and hence $\Gamma$ satisfies the property $Q_2.$

This completes the proof when $n\equiv 2\,(mod\,\,4).$

\NI Case 2: $n\equiv 0\,(mod\,\,4).$

Since $n\equiv 0\,(mod\,\,4)$ and $\Gamma\in \Gamma(2,\,1),$ $c$ is always odd, otherwise, the $b$-edges will not induce a Hamilton cycle. As $\alpha=2$ and $n\equiv 0\,(mod\,\,4),$ $a\equiv 2\,(mod\,\,4)$ as $gcd(a,\,n)=\alpha=2.$ By assumption $b$ is odd and hence $b\equiv 1\,or\,3\,(mod\,\,4);$ then $-b\equiv 3\,or \,1\,(mod\,\,4).$ For any two vertices $x$ and $y$ on $C_1,$ the path from $x$ to $y$ along $C_1$ has even length if and only if $x\equiv y\,(mod\,\,4).$ Hence the length of the section of $C_1$ from $b$ to $(-b)$ (containing the vertex $(b+a))$ is odd and the length of the section of $C_1$ from the vertex $(-b+a)$ to $b$ (containing the vertex $(b-a))$ is even and let this path be $R,$ (see Figure $4(b);$ note that the Figure $4(b)$ is for the case $n\equiv 2\,(mod\,\,4);$ the figure for the case $n\equiv 0\,(mod\,\,4)$ is similar). Now the edges of $R$ together with the edges $\{(b+a)b,\,(-b+a)a\}$ induce an even length path, say $R_1,$ from the vertex $b+a$ to $a.$ Again, the length of the section of $C_0$ from $a$ to $-2b+a$ (containing the vertex $2a$) is odd since it is of same length as the section of $C_1$ from $b$ to $-b$ (containing the vertex $(b+a))$. Now the edges of the section of $C_0$ from $a$ to $-2b+a$ together with the edge $a(-b+a)$ induce an even length path, say $R_2.$

As in Case 1, except the two $4$-cycles, namely, $(0,\,a,\,a+b,\,b)$ and $(-2b,\,-2b+a,\,-b+a,\,-b)$ in $\Gamma,$ each of the $4$-cycles formed by two consecutive vertices of $C_0$ and their corresponding vertices in $C_1$ constitute an alternating $4$-cycle with respect to $H_1$ and $H_2.$ Clearly, one of the two paths along $H_1,$ joining a pair of opposite vertices of these alternating $4$-cycles contains exactly one of the even length paths $R_1$ or $R_2$ and hence its length along $H_1$ is odd; hence if $(x,\,x+a,\,x+a+b,\,x+b)$ is an alternating $4$-cycle with $x,\,x+a\in V(C_0)$ and $x+a+b,\,x+b\in V(C_1),$ then $p_{H_1}(x,\,x+a+b)$ is odd.

We now show that in at least one of these alternating $4$-cycles, described above, $p_{H_2}(x,\,x+a+b)$ is odd. Let $S_0$ and $S_1$ be as defined in Case 1. Suppose the vertices of $C_0$ are alternately in $S_0$ and $S_1,$ then $c$ is even, but it is not the case. Therefore, there exist two consecutive vertices along $C_0$ which are in $S_0$ or $S_1.$ Thus there exists a pair of opposite vertices of an alternating $4$-cycle $(x,\,x+a,\,x+a+b,\,x+b),$ where $x\in V(C_0)$ and $x+a+b\in V(C_1)$ such that $p_{H_2}(x,\,x+a+b)$ is odd and hence $(x,\,x+a,\,x+a+b,\,x+b)$ is an odd alternating $4$-cycle. Thus $\Gamma$ satisfies property $Q_2,$ with respect to $H_1$ and $H_2.$

This completes the proof of the lemma.\hfill$\Box$

We use the following remarks in the proof of Lemma \ref{LEM1.directcirculant} given below.

\begin{lem}\label{lemma1ofbermond}\emph{\cite{JCTB.46.142-153}}
Let $\Gamma$ be a graph of $\Gamma(\alpha + 2,\beta).$ If the reduced graph $\Gamma^\prime$
admits a hamiltonian decomposition having the property $Q_1$ between $C_{\alpha-1}$  and
$C_0,$ then $\Gamma$ admits a Hamilton cycle decomposition having the property $Q_1$
between $C_{\alpha+1}$ and $C_0.$\hfill$\Box$
\end{lem}

\begin{rmk}\label{DC.remark1}
In the construction of two edge disjoint Hamilton cycles $H_1$ and $H_2$ of $\Gamma\in\Gamma(\alpha+2,\,\beta)$ from the two edge disjoint Hamilton cycles $H_1^\prime$ and $H_2^\prime,$ respectively, of $\Gamma^\prime\in\Gamma(\alpha,\,\beta),$ every $b$-edge of $H_i^\prime,\,i=1,\,2,$ connecting the vertices of $C_{\alpha-1}^\prime$ to $C_0^\prime$ is replaced by a path of odd length, see Lemma \ref{lemma1ofbermond}. Hence the parity of length of the path between any pair of vertices in $C_0,\,C_1,\,\ldots,\,C_{\alpha-1}$ along $H_i$ remains the same as the parity of length of the path between the respective vertices in $C_0^\prime,\,C_1^\prime,\,\ldots,\,C_{\alpha-1}^\prime$ of $H_i^\prime,$ where $C_i^\prime$ are the $a$-cycles of $\Gamma^\prime.$\hfill$\Box$
\end{rmk}

By $\overrightarrow{(a,\,b)}$ we denote a directed arc with tail at $a$ and head at $b.$
\begin{rmk}\label{DC.remark2}
Consider the graph $\Gamma\in\Gamma(3,\,2)$ and its reduced graph $\Gamma^\prime\in\Gamma(1,\,2).$ Let $\{H_1^\prime,\,H_2^\prime\}$ be the Hamilton cycle decomposition of $\Gamma^\prime$ guaranteed by Lemma \ref{LEM2.directcirculant}. Clearly, $H_1^\prime$ has two natural orientations one in the clockwise direction and the other in the anticlockwise direction and these two orientations induce orientations for the two $b$-edges of $H_1^\prime,$ namely, $0(-b)$ and $a(-b+a).$ We first fix the clockwise orientation of the Hamilton cycle $H_1^\prime$ of $\Gamma^\prime,$ and let the orientation of the edges $0(-b)$ and $a(-b+a)$ of $H_1^\prime$ give two arcs $\overrightarrow{(-b,\,0)}$ and $\overrightarrow{(-b+a,\,a)},$ see Figure 5. With respect to this orientation, we obtain a Hamilton cycle decomposition $\{H_1,\,H_2\}$ of $\Gamma,$ using the proof of Lemma \ref{lemma1ofbermond}, with property $Q_1,$ where the arcs $\overrightarrow{(-b,\,0)}$ and $\overrightarrow{(-b+a,\,a)}$ correpond to the $(m_1=2)$ edges $\{(\alpha-1,\,j_0)(0,\,j_0+c)\}$ and $\{(\alpha-1,\,j_1)(0,\,j_1+c)\},$ respectively, in the proof of Lemma \ref{lemma1ofbermond}.

\begin{center}
\includegraphics{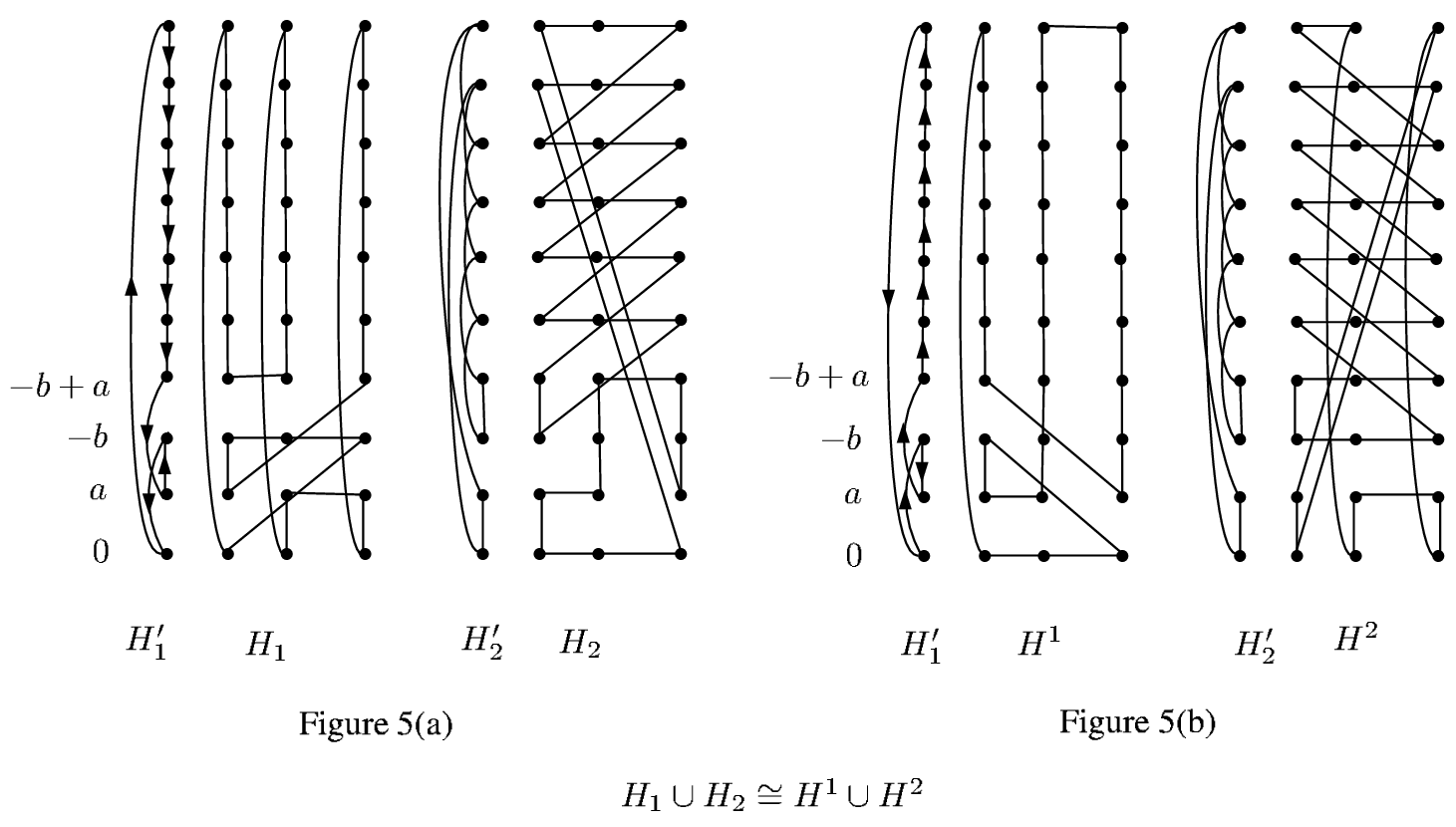}
\end{center}
{\small The Hamilton cycles $H_1$ and $H_2$ shown in Figure 5(a) arise out of the Hamilton cycles $H_1^\prime$ and $H_2^\prime$ with respect to the clockwise orientation of $H_1^\prime$ and, the Hamilton cycle decomposition $\{H^1,\,H^2\}$ shown in Figure 5(b) arises out of the Hamilton cycles $H_1^\prime$ and $H_2^\prime$ with respect to the anticlockwise orientation of $H_1^\prime,$ where $|\Gamma^\prime|=n=10,\,a=3,\,b=4.$}\\\vspace{-0.7cm}\begin{center}Figure 5\end{center}

If we fix the anticlockwise orientation for $H_1^\prime,$ we obtain another orientation of the edges $0(-b)$ and $a(-b+a)$ of $H_1^\prime$ and the resulting arcs are $\overrightarrow{(0,\,-b)}$ and $\overrightarrow{(a,\,-b+a)},$ see Figure $5(b).$ With respect to this orientation, we obtain a Hamilton cycle decomposition $\{H^1,\,H^2\}$ of $\Gamma,$ by proof of Lemma \ref{lemma1ofbermond}, with the property $Q_1,$ where the arcs $\overrightarrow{(0,\,-b)}$ and $\overrightarrow{(a,\,-b+a)}$ correspond to the $(m_1=2)$ edges $\{(\alpha-1,\,j_0)(0,\,j_0+c)\}$ and $\{(\alpha-1,\,j_1)(0,\,j_1+c)\},$ respectively, in the proof of Lemma \ref{lemma1ofbermond}, see Figure 5. Irrespective of the orientations of $H_1^\prime$ (that is, clockwise or anticlockwise), $H_1\cup H_2\cong H^1\cup H^2,$ that is, the two circulant graphs obtained from the union of the Hamilton cycles $H_1\cup H_2$ and $H^1\cup H^2$ in two different orientations 	of the edges of $H_1^\prime$ are isomorphic to each other. This isomorphism can be described by mapping $(i,\,j)$ to $(i,\,n-j),\,0\le i\ne \alpha-1,\,1\le j\le k_a-1,$ where the even integer $n$ is the number of vertices of $\Gamma^\prime,$ and the vertices $(i,\,0),0\le i\le \alpha-1$, are the fixed vertices of the isomorphism.\hfill$\Box$
\end{rmk}

\begin{lem}\label{LEM1.directcirculant}
Let $\Gamma\in\Gamma(\alpha+2,\,\beta)$ be a graph with even number of vertices and having the generating set $\{s,\,t\}$ where $s$ and $t$ are of different parity. If the reduced graph $\Gamma^\prime\in\Gamma(\alpha,\,\beta),$ of $\Gamma,$ has the properties $Q_1$ and $Q_2$ with respect to a Hamilton cycle decomposition $\{H_1^\prime,\,H_2^\prime\}$ of $\Gamma^\prime,$ then $\Gamma$ has a Hamilton cycle decomposition $\{H_1,\,H_2\}$ with the properties $Q_1$ and $Q_2.$
\end{lem}
\NI{\bf Proof.} Throughout the proof we use $\{s,\,t\}$ and $\{a,\,b\}$ as the generating sets for $\Gamma\in\Gamma(\alpha+2,\,\beta)$ and $\Gamma^\prime\in\Gamma(\alpha,\,\beta),$ respectively. Without loss of generality assume that throughout the proof $s$ is even and $t$ is odd. As $\Gamma$ is of even order and as $s$ is even and $t$ is odd, $\alpha$ is even and $\beta$ is odd, since $\alpha+2=gcd(m,\,s)$ and $\beta=gcd(m,\,t),$ where $m$ is the number of vertices of $\Gamma.$ Let $\{H_1^\prime,\,H_2^\prime\}$ be a Hamilton cycle decomposition of $\Gamma^\prime$ and $\{H_1,\,H_2\}$ be the corresponding Hamilton cycle decomposition of $\Gamma.$

In the construction of $H_1$ and $H_2$ from $H_1^\prime$ and $H_2^\prime$ (see \cite{JCTB.46.142-153}), $H_i,\,i=1,\,2,$ is obtained from $H_i^\prime$ by replacing each of the edges of $H_i^\prime$ between $C_{\alpha-1}^\prime$ to $C_0^\prime$ by an odd length path, whose origin is in $C_{\alpha-1}$ and terminus is in $C_0$ of $\Gamma$ and the internal vertices of the paths are in $C_\alpha$ and $C_{\alpha+1}.$

As $\{H_1^\prime,\,H_2^\prime\}$ has the property $Q_1,$ there exists a Hamilton cycle decomposition $\{H_1,\,H_2\}$ of $\Gamma\in\Gamma(\alpha+2,\,\beta)$ satisfying property $Q_1$ by Lemma \ref{lemma1ofbermond}. Next we prove that $\Gamma$ satisfies the property $Q_2$ with respect to $\{H_1,\,H_2\}.$ We prove this by induction on $\alpha.$

First we explain the idea behind the proof of this Lemma. In Claim 1 below, we prove the existence of a Hamilton cycle decomposition of $\Gamma\in\Gamma(3,\,2)$ satisfying $Q_2$ from a Hamilton cycle decomposition of $\Gamma^\prime\in\Gamma(1,\,2)$ with property $Q_2$ and in Claim 2 below, we prove the existence of a Hamilton cycle decomposition of $\Gamma\in\Gamma(4,\,1)$ satisfying $Q_2$ from a Hamilton cycle decomposition of $\Gamma^\prime\in\Gamma(2,\,1)$ with property $Q_2.$ In Claim 3 below, we obtain a Hamilton cycle decomposition of $\Gamma\in\Gamma(\alpha+2,\,\beta)$ with property $Q_2$ from a Hamilton cycle decomposition of $\Gamma^\prime\in\Gamma(\alpha,\,\beta)$ with property $Q_2.$ First we assume that $\{H_1^\prime,\,H_2^\prime\}$ is the Hamilton cycle decomposition of $\Gamma^\prime\in\Gamma(\alpha,\,\beta),$ where $\{\alpha,\,\beta\}=\{1,\,2\},$ as described in Lemma \ref{LEM2.directcirculant}.

\NI{\bf Claim 1.} For a Hamilton cycle decomposition $\{H_1^\prime,\,H_2^\prime\}$ of $\Gamma^\prime\in\Gamma(1,\,2)$ with property $Q_2,$ there is a Hamilton cycle decomposition $\{H_1,\,H_2\}$ of $\Gamma\in\Gamma(3,\,2)$ with property $Q_2.$

Let $|V(\Gamma^\prime)|=n;$ by hypothesis $\Gamma^\prime$ admits a Hamilton cycle decomposition $\{H_1^\prime,\,H_2^\prime\}$ with the property $Q_2.$ We consider two cases.

\NI{\bf Case 1.} $n\equiv 0\,(mod\,\,4).$

Let $\Gamma^\prime\in\Gamma(1,\,2)$ and let its corresponding graph in $\Gamma(3,\,2)$ be $\Gamma.$ From the definition of $\Gamma(\alpha,\,\beta),$ the labels of the vertices of the $i$th $a$-cycle are $(i,\,j),\,0\le i\le \alpha-1,\,0\le j\le k_a-1,$ where $k_a$ is the length of the $i$th $a$-cycle. We know that $\Gamma^\prime\in\Gamma(1,\,2)$ has only one $a$-cycle and hence its vertices are $(0,\,j),\,j=0,\,1,\,\ldots,\,n-1.$ Each vertex $(0,\,j)$ of $\Gamma^\prime$ gives rise to three vertices, namely, $(0,\,j),\,(1,\,j)$ and $(2,\,j),$ in $\Gamma$ and we call these three vertices of $\Gamma$ as the corresponding vertices of $(0,\,j)$ of $\Gamma^\prime$ and vice versa.

\begin{center}
\includegraphics[scale=0.9]{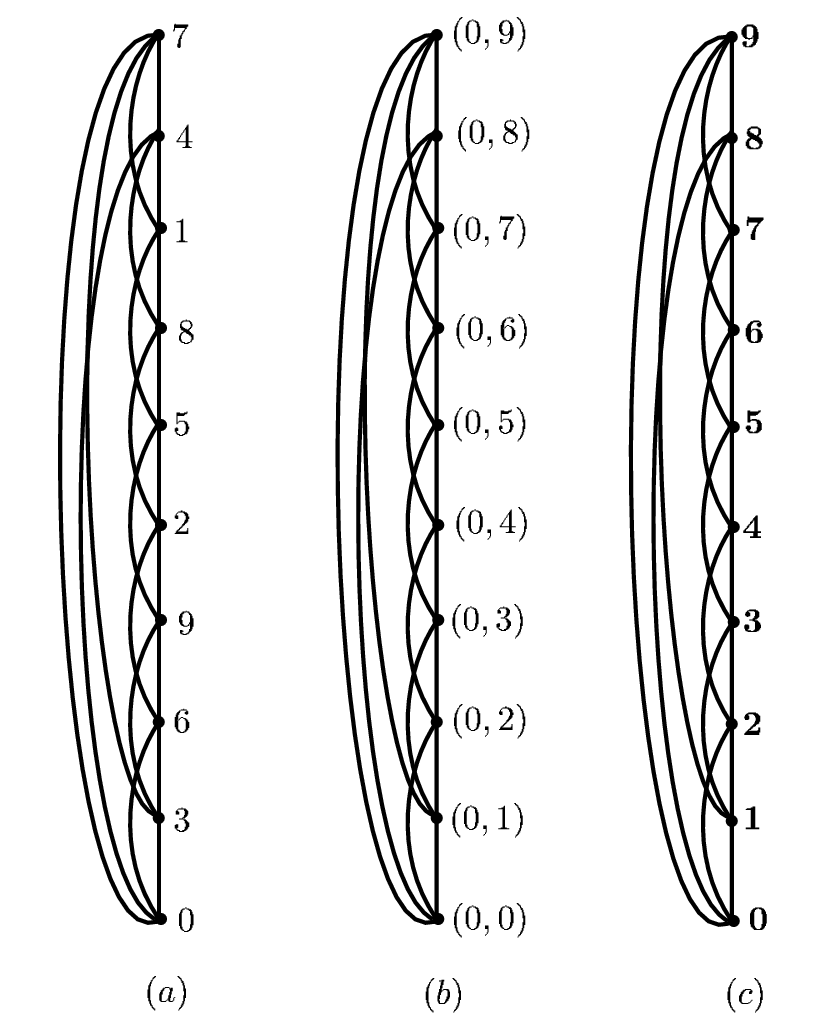}
\end{center}
{\small $(a)$ Labeling of the vertices of the circulant graph $\Gamma^\prime\in\Gamma(1,\,2)$ on $10$ vertices with the group elements, where the generating set is $\{3,\,4\}\subset \mathbb{Z}_{10}$ $(b).$ Labeling of the vertices of the graph $\Gamma^\prime$ with ordered pairs as described in \cite{JCTB.46.142-153}. $(c).$ Labeling of the vertices of the graph $\Gamma^\prime,$ with ${\bf j}$ if the vertex has the label $(0,\,j)$ in $(b).$}\vspace{-.4cm}\begin{center}Figure 6
\end{center}

For our convenience we relabel the vertices $(0,\,j),\,0\le j\le n-1,$ of $\Gamma^\prime\in\Gamma(1,\,2)$ as ${\bf j},$ see Figure 6(c); we call it as {\it new labelling of $\Gamma^\prime$}(throughout this lemma, the new labelling of $\Gamma^\prime$ is denoted by bold face letters);  in fact, in the proof of this lemma each vertex of $\Gamma^\prime$ will have three different labels, namely, $(0,\,j),\,{\bf j}$ and the other one the group element and according to our convenience and circumstances we use one of these labels. But, for the vertices of $\Gamma\in\Gamma(3,\,2)$ we use the unique label $(i,\,j)$ as in \cite{JCTB.46.142-153}. In the graph $\Gamma^\prime,$ let the vertex $-b$ be denoted by the label $(0,\,r)$ and hence in our new labeling it is denoted by ${\bf r},$ see Figure 6; for our convenience we write $-b={\bf r}.$ Clearly, ${\bf r+1}\,(=-b+a)$ is the immediate next vertex of ${\bf r}\,(=-b)$ in $C_0^\prime$ of $\Gamma^\prime$ in the new labeling. Corresponding to the vertices ${\bf r}$ and ${\bf r+1}$ of $\Gamma^\prime,$ there are two rows, each having three vertices, in $\Gamma;$ the vertices of these rows are denoted by, $(0,\,r),\,(1,\,r),\,(2,r)$ and $(0,\,r+1),\,(1,\,r+1),\,(2,r+1),$ respectively, see Figure 7.

\begin{center}
\includegraphics{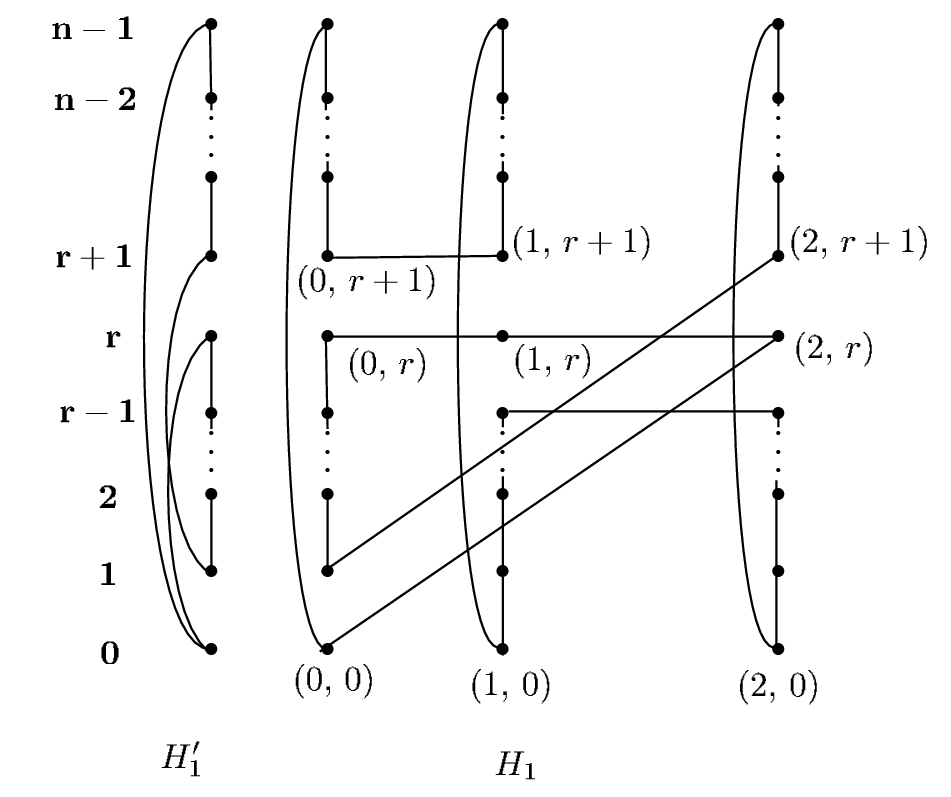}
\end{center}
{\small Construction of the Hamilton cycle $H_1$ of $\Gamma\in\Gamma(3,\,2)$ from the Hamilton cycle $H_1^\prime$ of $\Gamma^\prime\in\Gamma(1,\,2),$ as described in \cite{JCTB.46.142-153}.}\begin{center}
Figure 7
\end{center}

In the construction of $\Gamma\in\Gamma(3,\,2)$ from $\Gamma^\prime\in\Gamma(1,\,2),$ the two edges $(0,\,r)(1,\,r)$ and $(0,\,r+1)(1,\,r+1)$ are in $H_1$ and the two edges $(0,\,r)(0,\,r+1)$ and $(1,\,r)(1,\,r+1)$ are in $H_2,$ where $\{H_1,\,H_2\}$ is the Hamilton cycle decomposition of $\Gamma,$ corresponding to the Hamilton cycle decomposition $\{H_1^\prime,\,H_2^\prime\}$ of $\Gamma^\prime,$ as described in \cite{JCTB.46.142-153}, see Figures 7 and 9. We claim that $\Gamma$ satisfies the property $Q_2$ with respect to the Hamilton cycle decomposition $\{H_1,\,H_2\},$ where the required $4$-cycle of $\Gamma$ is $(0,\,r)\,(0,\,r+1)\,(1,\,r+1)\,(1,\,r).$

First we show that $p_{H_1}((0,\,r),\,(1,\,r+1))$ is odd. Clearly, $C_0^\prime=({\bf 0},\,{\bf 1},\,{\bf 2},\,\ldots,\,{\bf n-1},\,{\bf 0})$ is a Hamilton cycle in $\Gamma^\prime,$ with respect to our new labeling.

The Hamilton cycle $H_1$ of $\Gamma$ contains the four vertices $(0,\,r),\,(0,\,0),\,(0,\,r+1)$ and $(1,\,r+1)$ in the clockwise order as shown in the Figure 7 ( the order is guaranteed by the corresponding Hamilton cycle $H_1^\prime$ of $\Gamma^\prime$), where $(0,\,r+1)(1,\,r+1)$ is an edge of $H_1.$ As $C_0^\prime$ being an even cycle, it can be thought of as a bipartite graph with bipartition $X=\{{\bf 0},\,{\bf 2},\,{\bf 4},\,\ldots,\,{\bf n-2}\}$ and $Y=\{{\bf 1},\,{\bf 3},\,\ldots,\,{\bf n-1}\}.$ Clearly, the vertex $-b\,\big(={\bf r},$ the $({\bf r+1})$th vertex along $C_0^\prime,$ the group element $ra\,(mod\,\,n);$ we do not differentiate the labels $-b,\,{\bf r}$ and $ra$ and we denote the vertex $-b$ by $-b={\bf r}=ra$ in the three labellings of the vertices of $\Gamma^\prime\big)$ is in $X,$ because $b$ is even implies $-b$ is even; as $-b=ra$ and $a$ is odd, ${\bf r}$ is even. Further, $0\in X$ and $-b\in X$ implies $p_{C_0^\prime}({\bf 0},\,{\bf r})=p_{C_0^\prime}(0,\,-b)$ is even. Now consider the Hamilton cycle decomposition $\{H_1^\prime,\,H_2^\prime\}$ of $\Gamma^\prime$ as in Lemma \ref{LEM2.directcirculant}. The Hamilton cycle $H_1^\prime$ of $\Gamma^\prime$ is obtained from $C_0^\prime$ by deleting two edges $0a\,(={\bf 01})$ and $(-b)(-b+a)(={\bf r(r+1)})$ of $C_0^\prime$ and adding two edges $0(-b)(={\bf 0r})$ and $a(-b+a)(={\bf 1}({\bf r+1})),$ see Figure 8.

\begin{center}
\includegraphics[scale=0.9]{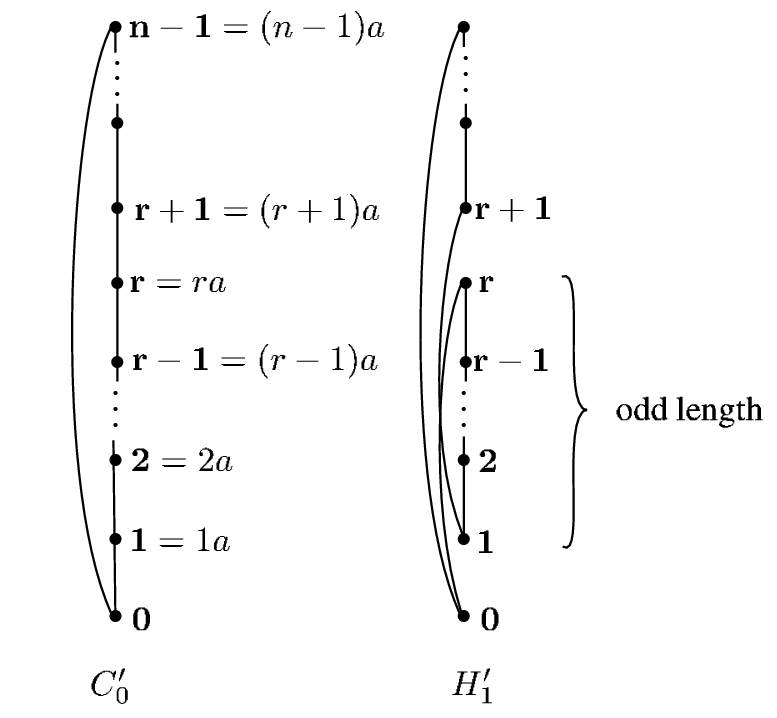}
\end{center}
{\small Length of the $({\bf 1},\,{\bf r})$-section along $H_1^\prime$ is odd, as ${\bf r}$ (an even integer) is the length from ${\bf 0}$ to ${\bf r}$ along $C_0^\prime;$ this implies that the  $({\bf 0},\,{\bf r+1})$-section along $H_1^\prime$ is of odd length.}
\begin{center}Figure 8
\end{center}

Next we prove that $p_{H_1^\prime}({\bf 0},\,{\bf r+1})\,(=p_{H_1^\prime}(0,\,-b+a))$ is odd. From the last paragraph, $p_{C_0^\prime}({\bf 0},\,{\bf r})$ is even and hence $p_{C_0^\prime}({\bf 1},\,{\bf r})$ is odd, see Figure 8; the path from ${\bf 1}$ to ${\bf r}$ along $C_0^\prime$ is also in $H_1^\prime$ and hence $p_{H_1^\prime}({\bf 1},\,{\bf r})$ is odd; consequently, as $n$ is even,\begin{eqnarray} p_{H_1^\prime}({\bf 0},\,{\bf r+1})\,is\,odd,\,\label{directcirculant.equation1}\end{eqnarray} (where we consider the $({\bf 0},\,{\bf r+1})$-section of $H_1^\prime$ not containing the vertex ${\bf r}$), see Figure 8.

Next we prove that $p_{H_1}((0,\,r),\,(1,\,r+1))$ is odd. We divide the $((0,\,r),\,(1,\,r+1))$-section of $H_1$ (containing the vertex $(1,\,r))$ into three subsections and we show that each one of them is of odd length to conclude $p_{H_1}((0,\,r),\,(1,\,r+1))$ is odd; the subsections are $((0,\,r),\,(0,\,0)),$ $((0,\,0),\,(0,\,r+1))$ and $((0,\,r+1),\,(1,\,r+1)),$ see Figure 7. From the construction of $H_1$ from $H_1^\prime,$ it is clear that as $((0,\,r),\,(0,\,0))$-section and $((0,\,r+1),\,(1,\,r+1))$-section of $H_1$ are of length 3 and length 1, respectively, it is enough to show that $((0,\,0),\,(0,\,r+1))$-subsection is of odd length. But it is easy to see that $((0,\,0),\,(0,\,r+1))$-section of $H_1$ is identical with the $({\bf 0},\,{\bf r+1})$-section of the Hamilton cycle $H_1^\prime,$ see Figure 7. The $({\bf 0},\,{\bf r+1})$-section of $H_1^\prime$ is already proved to be of odd length, by (\ref{directcirculant.equation1}). Thus $p_{H_1}((0,\,0),\,(0,\,r+1))$ is odd.

Next we prove that $p_{H_2}((0,\,r),\,(1,\,r+1))$ is odd. As above, we divide the  $((0,\,r),\,(1,\,r+1))$-section of $H_2$ (not containing the vertex $(0,\,r+1)$) into four subsections, namely, $((0,\,r),\,(0,\,0)),$\,((0,\,0),\,(0,\,1)),$\,((0,\,1),\,(0,\,r-1))$ and $((0,\,r-1),\,(1,\,r+1))$-section in the cyclic order are guaranteed by $H_2^\prime,$ see \cite{JCTB.46.142-153} and Figure 9. From the construction of $H_2,$ $(0,\,0)(0,\,1)$ is an edge and the $((0,\,r-1),\,(1,\,r+1))$-section is a path of length 3, namely, $(0,\,r-1)\,(1,\,r-1)\,(1,\,r)\,(1,\,r+1),$ see Figure 9. Hence we show that the lengths of the other two sections are of different parity. This is achieved by finding the lengths of the corresponding sections in $H_2^\prime$ of $\Gamma^\prime.$ The sections in $H_2^\prime$ corresponding to the $((0,\,r),\,(0,\,0))$-section and $((0,\,1),\,(0,\,r-1))$-section of $H_2$ in $\Gamma$ are $({\bf r},\,{\bf 0})$-section and $({\bf 1},\,{\bf r-1})$-section, respectively, in the new labeling of $\Gamma^\prime.$

\begin{center}
\includegraphics[scale=0.8]{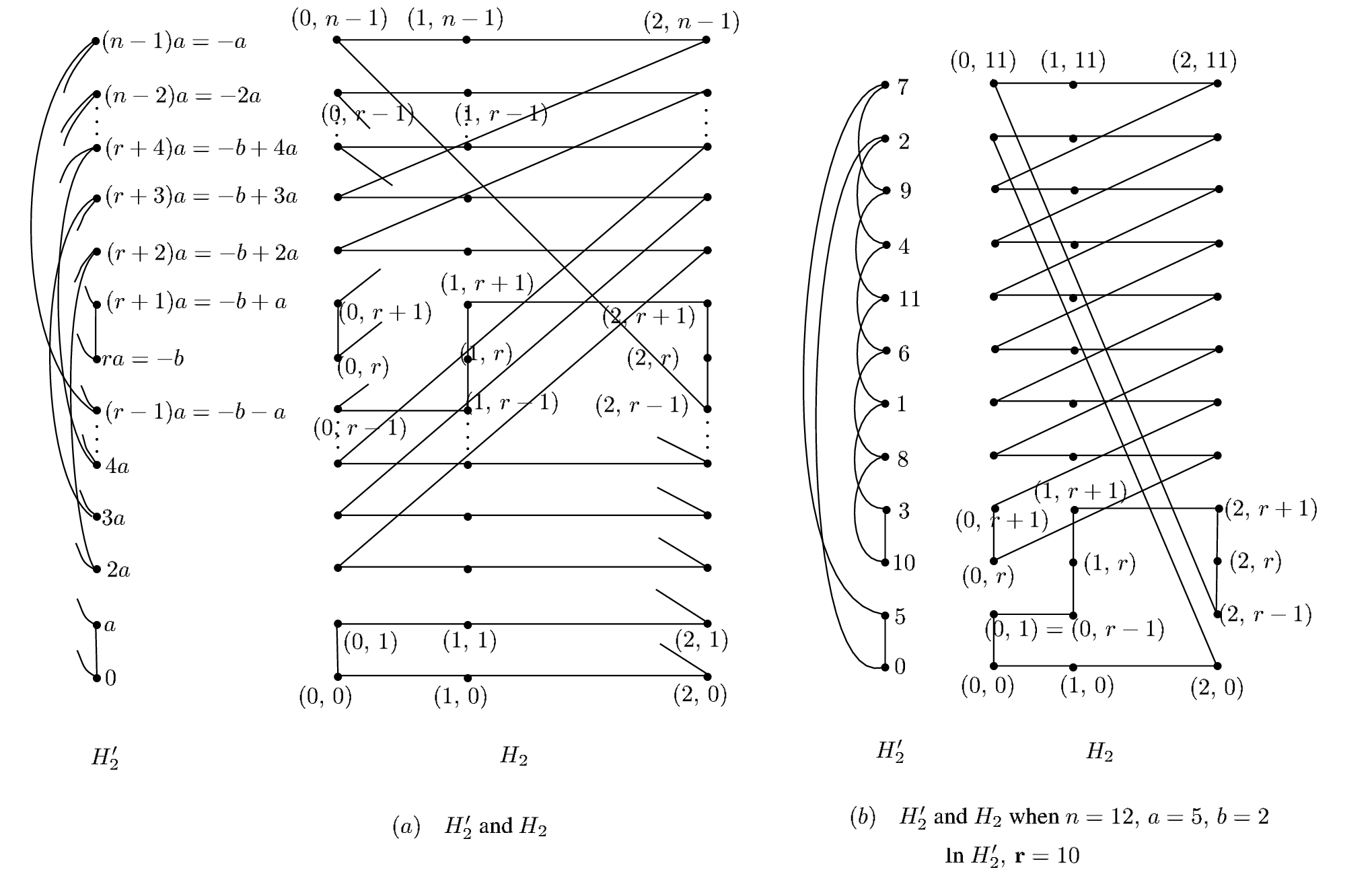}
{\small $H_2^\prime$ in $\Gamma^\prime$ and $H_2$ in $\Gamma$}\end{center}
\vspace{-.5cm}
\begin{center}Figure 9
\end{center}

We shall show that $p_{H_2^\prime}({\bf r},\,{\bf 0})$ is odd and $p_{H_2^\prime}({\bf 1},\,{\bf r-1})$ is even. We prove that $p_{H_2^\prime}({\bf r},\,{\bf 0})$ is odd by observing that the $({\bf 0},\,{\bf r})$-section (not containing the edge ${\bf 0r}$) of one of the $b$-cycles, namely $0\,b\,(2b)\,\ldots\,(-b)\,0\,\,\big(=({\bf 0}\,\ldots {\bf r}\,{\bf 0})\text{ in the new labeling,}$ not containing the edge ${\bf 0r}\big)$ of $\Gamma^\prime$ is a section of $H_2^\prime,$ and it is of odd length as the $({\bf 0},\,{\bf r})$-section of $H_2^\prime$ together with the edge ${\bf 0r}$ is a $b$-cycle, which is of even length in $\Gamma^\prime,$ that is, $p_{H_2^\prime}({\bf r},\,{\bf 0})$ is odd.

\begin{center}
\includegraphics{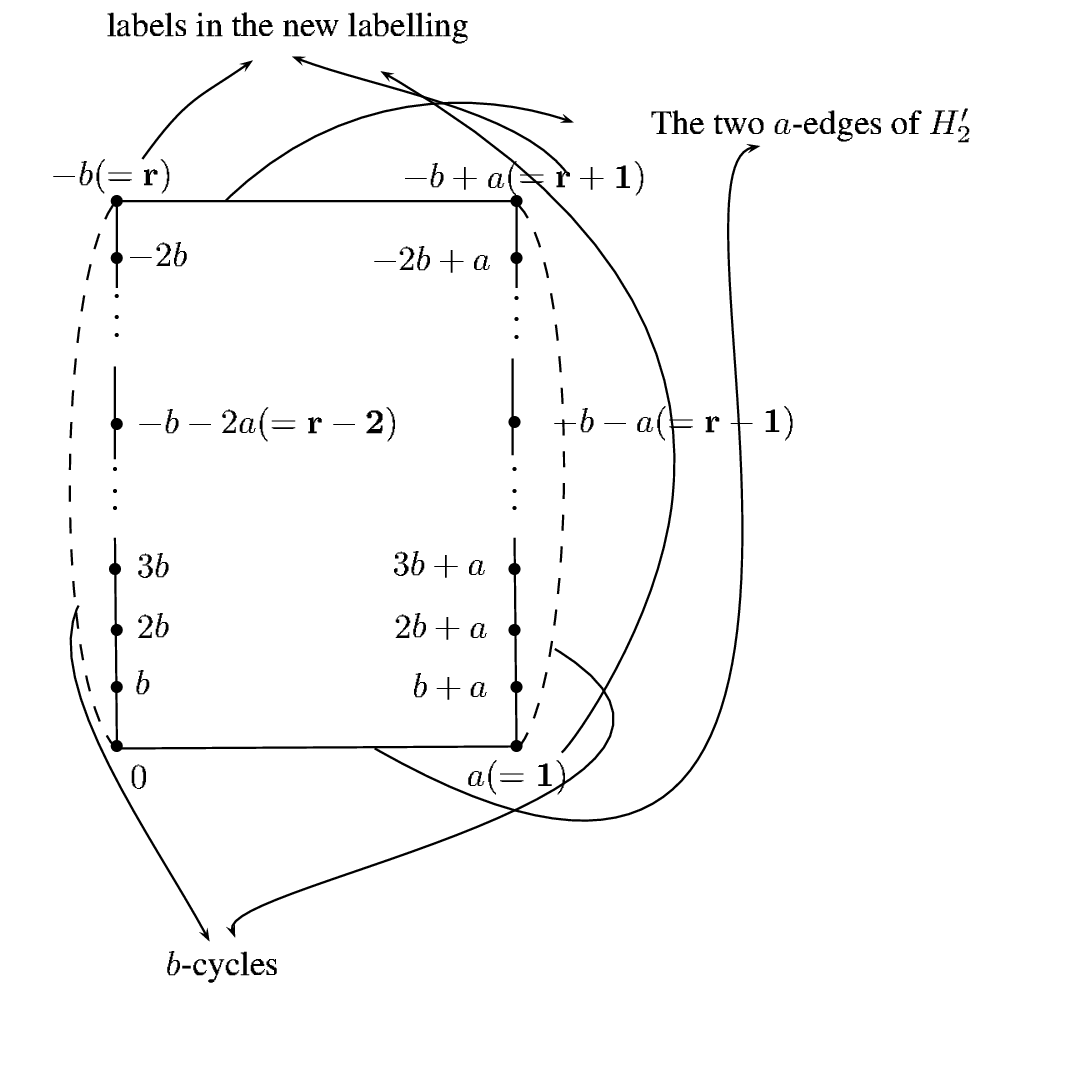}

{\small A redrawing of $H_2^\prime$}\end{center}
\begin{center}Figure 10\end{center}

Next we prove that $p_{H_2^\prime}({\bf 1},\,{\bf r-1})$ is even. Now consider the $({\bf 1},\,{\bf r-1})$-section of $H_2^\prime$ containing only $b$-edges. As we move along $b$-edges of the $({\bf 1},\,{\bf r-1})$-section of $H_2^\prime,$ the alternate labels of the vertices are congruent to $1\,(mod\,\,4)$ or $3\,(mod\,\,4).$ If two nonconsecutive vertices of $H_2^\prime$ are both congruent to $1\text{ or }3\,(mod\,\,4),$ then their distance along $H_2^\prime$ is even. Therefore, as $a\equiv 1\text{ or }3\,(mod\,\,4)$ implies $-b-a\equiv 1\text{ or }3\,(mod\,\,4),$ $p_{H_2^\prime}({\bf 1},{\bf r-1})$ is even.

As observed earlier, the length of the $((0,\,r),\,(1,\,r+1))$-section of $H_2$ is sum of the lengths of the sections $((0,\,r),\,(0,\,0)),$\,$((0,\,0),\,(0,\,1)),$\,$((0,\,1),\,(0,\,r-1))$ and a path of length $3$ from $(0,\,r-1)$ to $(1,\,r+1)$ in $H_2$ of $\Gamma.$ As $p_{H_2^\prime}({\bf r},\,{\bf 0})$ is odd, $p_{H_2}((0,\,r),\,(0,\,0))$ is odd by Remark \ref{DC.remark1}, the $((0,\,0),\,(0,\,1))$-section is an edge, the length of $((0,\,1),\,(0,\,r-1))$-section of $H_2$ is even, as $p_{H_2^\prime}({\bf 1},\,{\bf r-1})$ is even and by Remark \ref{DC.remark1} and the last section is a $P_4=(0,\,r-1)(1,\,r-1)(1,\,r)(1,\,r+1).$ Thus in $\Gamma,$ $p_{H_2}((0,\,r),\,(1,\,r+1))$ is odd.

\NI{\bf Case 2.} $n\equiv 2\,(mod\,\,4).$

Let $\{H_1^\prime,\,H_2^\prime\}$ be the Hamilton cycle decomposition of $\Gamma^\prime$ with property $Q_2$ obtained in Lemma \ref{LEM2.directcirculant}. As pointed out in Remark \ref{DC.remark2}, the Hamilton cycle $H_1^\prime$ has two natural orientations and with respect to each of these two orientations there is a Hamilton cycle decomposition of $\Gamma\in\Gamma(3,\,2).$ We show that in at least one of these two orientations of $H_1^\prime,$ the corresponding Hamilton cycle decomposition of $\Gamma$ has the property $Q_2$ (Note that both these two Hamilton cycle decompositions satisfy the property $Q_1$ as in Remark \ref{DC.remark2}). Consider the new labeling of $\Gamma^\prime$ described in Case 1 above.

\NI{\bf Subcase 2.1.} $p_{H_2^\prime}({\bf 1},\,{\bf r-1})$ is odd.

We know that $H_1^\prime$ has exactly two $b$-edges. In this case we consider the clockwise orientation of $H_1^\prime.$ Consider the $4$-cycle $C=((0,\,r)(0,\,r+1)(1,\,r+1)(1,\,r))$ in $\Gamma,$ see Figure 9 (the figure for the case $n\equiv 2\,(mod\,\,4)$ also resembles the same as in the case $n\equiv 0\,(mod\,\,4)$ and note that ${\bf r}$ is even).


To prove $p_{H_1}((0,\,r),\,(1,\,r+1))$ is odd, we divide the $((0,\,r),\,(1,\,r+1))$-section of $H_1,$ (containing the vertex $(1,\,r)$) into three subsections and we show that each of them is of odd length. The subsections are $((0,\,r),\,(0,\,0)),$ $((0,\,0),\,(0,\,r+1))$ and $((0,\,r+1),\,(1,\,r+1)),$ see Figure 7 (the figure for the case $n\equiv 2\,(mod\,\,4)$ also resembles the same as in the case $n\equiv 0\,(mod\,\,4)$). From the construction of $H_1$ from $H_1^\prime,$ it is clear that $((0,\,r),\,(0,\,0))$-section of $H_1$ is of length 3, see Figure 7, and $(0,\,r+1)(1,\,r+1)$ is an edge, hence it is enough to show that $((0,\,0),\,(0,\,r+1))$-section (containing the vertex $(0,\,n-1)$) is of odd length. But it is easy to see (in fact, it is the same proof as in the case $n\equiv 0\,(mod\,\,4)$) that this section of $H_1$ is the same as the $({\bf 0},\,{\bf r+1})$-section of the Hamilton cycle $H_1^\prime$ of $\Gamma^\prime,$ which is of odd length and hence $p_{H_1}((0,\,0),\,(0,\,r+1))$ is odd. This completes the proof that $p_{H_1}((0,\,r),\,(1,\,r+1))$ is odd.

Now we prove that $p_{H_2}((0,\,r)(1,\,r+1))$ is odd. For that we divide the $((0,\,r),\,(1,\,r+1))$-section of $H_2$ (not containing the vertex $(0,\,r+1)$) into four subsections, namely, $((0,\,r),\,(0,\,0)),\,((0,\,0),\,(0,\,1)),\,((0,\,1),\,(0,\,r-1))$ and $((0,\,r-1),\,(1,\,r+1)),$ see Figure 9. It is clear that the $((0,\,0),\,(0,\,1))$-section is an edge and the $((0,\,r-1),\,(1,\,r+1))$-section is a path of length $3.$ We shall show that the lengths of the subsections $((0,\,r),\,(0,\,0))$ and $((0,\,1),\,(0,\,r-1))$ are of different parity.

First we prove that the length of the subsection $((0,\,r),\,(0,\,0))$ of $H_2$ is even. Clearly, the section corresponding to $((0,\,r),\,(0,\,0))$ of $H_2$ is $({\bf r},\,{\bf 0})$-section in $H_2^\prime$ of $\Gamma^\prime,$ in the new labeling. Observe that the $({\bf 0},\,{\bf r})$-section of one of the $b$-cycles, not containing the edge ${\bf 0r},$ of $\Gamma^\prime$ is a section of $H_2^\prime,$ and it is of even length as the $({\bf 0},\,{\bf r})$-section of $H_2^\prime$ together with the edge ${\bf 0r}$ is a $b$-cycle of odd length in $\Gamma^\prime$ (as $n\equiv 2\,(mod\,\,4)$ and $\beta=2$), that is, $p_{H_2^\prime}({\bf r},\,{\bf 0})$ is even. Hence by Remark \ref{DC.remark1}, $p_{H_2}((0,r),\,(0,\,0))$ is even.

Next consider the $((0,\,1),\,(0,\,r-1))$-section (containing the vertex $(1,\,1)$) of $H_2;$ the corresponding $({\bf 1},\,{\bf r-1})$-section in $H_2^\prime$ contains only $b$-edges. As $p_{H_2^\prime}({\bf 1},\,{\bf r-1})$ is odd, by assumption, and also the section $({\bf 1},\,{\bf r-1})$ of $H_2^\prime$ containing only $b$-edges, $p_{H_2}((0,\,1),\,(0,\,r-1))$ is odd, by Remark \ref{DC.remark1}. This proves $p_{H_2}((0,\,r),\,(1,\,r+1))$ is odd.

\NI{\bf Subcase 2.2.} $p_{H_2^\prime}({\bf 1},\,{\bf r-1})$ is even.

In this subcase, we consider the anticlockwise orientation of the edges of $H_1^\prime.$ Then $C=((0,\,0)(0,\,1)(1,\,1)(1,\,0))$ will be proved to be an odd alternating $4$-cycle, see Figure 11.

First we prove that $p_{H_1}((0,\,0)(1,\,1))$ is odd. For this, we divide the $((0,\,0),\,(1,\,1))$-section of $H_1$ (containing the vertex $(1,\,0)$) into three subsections, namely, $((0,\,0),\,(0,\,r)),\,((0,\,r),\,(0,\,1))$ and $((0,\,1),\,(1,\,1)).$ Clearly, $((0,\,0),\,(0,\,r))$-section is the path $(0,\,0)\,(1,\,0)\,(2,\,0)\,(0,\,r)$ of length 3 and $((0,\,1),\,(1,\,1))$-section is an edge and hence it is enough to prove that $((0,\,r),\,(0,\,1))$-section of $H_1$ is of odd length.

In $H_1^\prime$ of $\Gamma^\prime,$ $b$ is even implies $-b$ is even as $n$ is even. As $-b=ra$ and $a$ is odd implies ${\bf r}$ is even. It is an easy observation that the vertices with even numbered labels in the new labeling of $H_1^\prime$ are at odd distance from the vertex $a={\bf 1},$ the new label of the vertex $a,$ since the even numbered vertices are $2a={\bf 2},\,4a={\bf 4},$ etc. Hence $p_{H_1^\prime}({\bf 1},\,{\bf r})$ is odd as ${\bf r}$ is even. As the $((0,\,r),\,(0,\,1))$-section of $H_1$ has the same length as the $({\bf r},\,{\bf 1})$-section of $H_1^\prime,$ $p_{H_1}((0,\,r),\,(0,\,1))$ is odd. 

Finally, we show that $p_{H_2}((0,\,0)(1,\,1))$ is odd. To prove this, we divide the $((0,\,0),\,(1,\,1))$-section of $H_2$ (containing the vertex $(0,\,1)$) into three subsections, namely, $((0,\,0),\,(0,\,1)),\,((0,\,1),\,(0,\,r-1))$ and $((0,\,r-1),\,(1,\,1)),$ see Figure 11. It is clear that the  $((0,\,0),\,(0,\,1))$-section is an edge and the $((0,\,r-1),\,(1,\,1))$-section is the path $(0,\,r-1)\,(2,\,n-1)\,(2,\,0)\,(2,\,1)\,(1,\,1)$ on $5$ vertices, see Figure 11. Hence it is enough to show that the $((0,\,1),\,(0,\,r-1))$-section is of even length. As $p_{H_2^\prime}({\bf 1},\,{\bf r-1})$ is even, by assumption in this subcase, $({\bf 1},\,{\bf r-1})$-section of $H_2^\prime,$ containing only $b$-edges, is of even length and hence $p_{H_2}((0,\,1),\,(0,\,r-1))$ is even, by Remark \ref{DC.remark1}. Thus $p_{H_2}((0,\,0),\,(0,\,1))$ is odd. Hence the Hamilton cycle decomposition $\{H_1,\,H_2\}$ of $\Gamma$ satisfies the property $Q_2.$

\begin{center}
\includegraphics{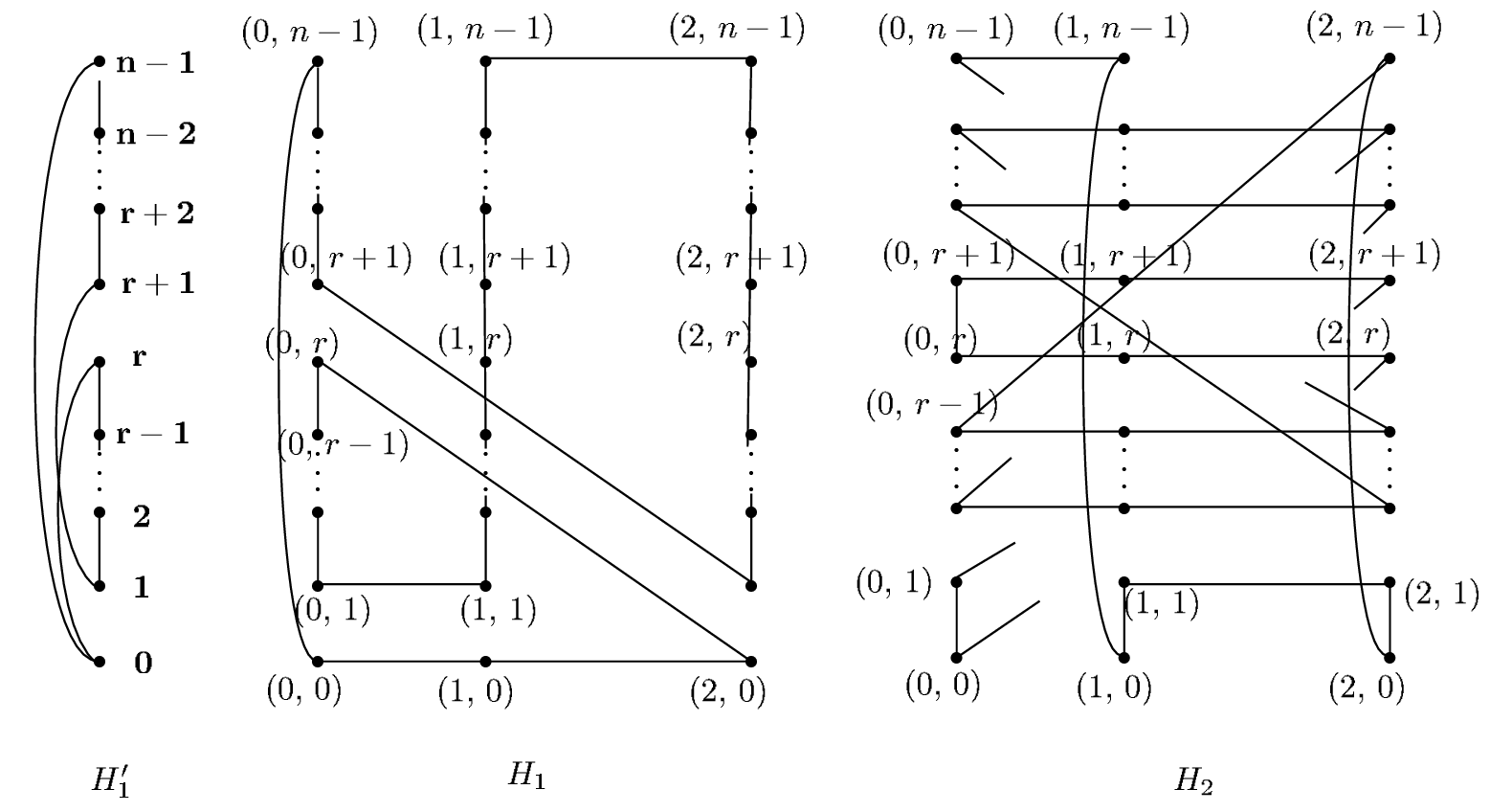}
\end{center}
{\small The Hamilton cycles $H_1$ and $H_2$ of $\Gamma\in\Gamma(3,\,2)$ corresponding to the anticlockwise orientation of $H_1^\prime$ of $\Gamma^\prime\in\Gamma(1,\,2).$}\begin{center}Figure 11
\end{center}

\NI{\bf Claim 2.} For any Hamilton cycle decomposition $\{H_1^\prime,\,H_2^\prime\}$ of $\Gamma^\prime\in\Gamma(2,\,1)$ with the property $Q_2,$ there is a Hamilton cycle decomposition $\{H_1,\,H_2\}$ of $\Gamma\in\Gamma(4,\,1)$ with the property $Q_2.$ 

Consider the graph $\Gamma^\prime\in\Gamma(2,\,1).$ There are two $a$-cycles $C_0^\prime$ and $C_1^\prime;$ hence there are two parallel matchings between $C_0^\prime$ and $C_1^\prime.$ We denote the matching from $C_0^\prime$ to $C_1^\prime$ as $M_0$ and the other matching from $C_1^\prime$ to $C_0^\prime$ as $M_1;$ the edges of $M_0$ are $(0,\,j)\,(1,\,j)$ and the edges of $M_1$ are $(1,\,j)\,(0,\,j+c),$ $0\le j\le \frac{n}{2}-1,$ for some $c.$

Clearly, from the construction of $H_1^\prime$ and $H_2^\prime$ of $\Gamma^\prime,$ as in \cite{JCTB.46.142-153}, every alternating $4$-cycle having two consecutive vertices along $C_0^\prime$ and two consecutive vertices along $C_1^\prime$ uses exactly two $a$-edges and two $b$-edges and both these two $b$-edges are in $M_0$ or in $M_1.$ Lemma \ref{LEM2.directcirculant} guarantees the existence of an odd alternating $4$-cycle with respect to the Hamilton cycle decomposition $\{H_1^\prime,\,H_2^\prime\}$ of $\Gamma^\prime$ and this $4$-cycle contains two $b$-edges of $M_0.$ In the construction of $H_i,\,i=1,\,2,$ of $\Gamma$ (from $H_i^\prime,\,i=1,\,2,$ of $\Gamma^\prime$) the $a$-edges and $M_0$-edges (that is, $b$-edges that belong to the matching $M_0$) of $\Gamma^\prime$ are retained as it is in the transformation of $\Gamma^\prime$ to $\Gamma.$ While obtaining $H_i$ from $H_i^\prime,$ each of the $b$-edges of $M_1$ in $H_i^\prime$ are replaced by a path of odd length whose internal vertices are in $C_2$ and $C_3.$ Hence an odd alternating $4$-cycle that exists in $\Gamma^\prime$ becomes an odd alternating $4$-cycle with respect to the Hamilton cycle decomposition $\{H_1,\,H_2\}$ of $\Gamma.$

\NI{\bf Claim 3.}  For any Hamilton cycle decomposition $\{H_1^\prime,\,H_2^\prime\}$ of $\Gamma^\prime\in\Gamma(\alpha,\,\beta)$ with property $Q_2,$ there is a Hamilton cycle decomposition $\{H_1,\,H_2\}$ of $\Gamma\in\Gamma(\alpha+2,\,\beta)$ with property $Q_2.$ 

In Claims 1 and 2, we proved that the existence of the property $Q_2$ in $\Gamma^\prime\in \Gamma(1,\,2)$ or $\Gamma(2,\,1)$ implies the existence of the property $Q_2$ in $\Gamma\in\Gamma(3,\,2)$ or $\Gamma(4,\,1),$ respectively. From the Hamilton cycle decomposition and construction of the odd alternating $4$-cycle of $\Gamma$ from $\Gamma^\prime,$ we have shown that the odd alternating $4$-cycle lies between $C_0$ and $C_1,$ that is having two vertices in $C_0$ and two vertices in $C_1.$ Hence, as in Claim 2, the odd alternating $4$-cycle in the Hamilton cycle decomposition $\{H_1^\prime,\,H_2^\prime\}$ of $\Gamma^\prime\in \Gamma(\alpha,\,\beta)$ is also the odd alternating $4$-cycle in the Hamilton cycle decomposition $\{H_1,\,H_2\}$ of $\Gamma\in\Gamma(\alpha+2,\,\beta);$ this completes the proof of existence of the property $Q_2$ with respect to the Hamilton cycle decomposition $\{H_1,\,H_2\}$ of $\Gamma.$

This completes the proof of the lemma.\hfill$\Box$

\begin{rmk}\label{DC.remark3}
Recall that our graphs of $\Gamma(\alpha,\,\beta)$ always have even order by assumption. As every member of $\Gamma(\alpha,\,\beta)$ is isomorphic to a member of $\Gamma(\beta,\,\alpha)$ and vice-versa, to each member of $\Gamma\in\Gamma(\alpha,\,\beta),$ there corresponds a graph $\Gamma_1\,(\cong\Gamma)\in\Gamma(\beta,\,\alpha).$ If $\alpha$ is odd, then $k_a,$ the length of each of the $a$-cycles of $\Gamma\in\Gamma(\alpha,\,2),$ must be even. Let $\Gamma\in\Gamma(\alpha,\,2)$ and let $\Gamma_1\in\Gamma(2,\,\alpha)$ be the image of $\Gamma$ under an isomorphism. Clearly, the $a$-edges and $b$-edges of $\Gamma$ are the $b$-edges and $a$-edges of $\Gamma_1,$ respectively, and vice versa. As $k_a$ is even, each $a$-cycle $C_i,\,0\le i\le \alpha-1,$ in $\Gamma$ has two $1$-factors, say $F_1^i$ and $F_2^i$ in the subgraph induced by $C_i;$ that is, for $i\in \{0,\,1,\,\ldots,\, \alpha-1\},$ $F_1^i=\{(i,\,0)(i,\,1),\,(i,\,2)(i,\,3),\,\ldots,\,(i,\,n-2)(i,\,n-1)\}$ and $F_2^i=\{(i,\,1)(i,\,2),\,(i,\,3)(i,\,4),\,\ldots,\,(i,\,n-1)(i,\,0)\}.$ The edges of the $1$-factor $F_1^i$ of $\Gamma$ become the $M_0$-edges of $\Gamma_1\in\Gamma(2,\,\alpha)$ and the edges of the $1$-factor $F_2^i$ of $\Gamma$ become the $M_1$-edges of $\Gamma_1.$ In the proof of Lemma \ref{LEM1.directcirculant}, we constructed an odd alternating $4$-cycle with respect to the Hamilton cycle decomposition of $\Gamma,$ so that the two $a$-edges use the edges of $F_1^i.$ Hence, if $\Gamma\in\Gamma(\alpha,\,2)$ has an odd alternating $4$-cycle, then the image of the $4$-cycle under the isomorphism between the corresponding graphs of $\Gamma(\alpha,\,2)$ and $\Gamma(2,\,\alpha)$ gives the odd alternating $4$-cycle in $\Gamma_1\in\Gamma(2,\,\alpha)$ using two $M_0$-edges of $\Gamma_1.$
\end{rmk}

\begin{rmk}\label{DC.remark4}
By saying, for $\Gamma\in\Gamma(\alpha,\,\beta)$ obtain the reduced graph of $\Gamma$ in $\Gamma(\alpha_1,\,\beta_1),\,\alpha_1\le\alpha,\,\beta_1\le\beta,$ we mean the following: The successive reduced graphs of $\Gamma\in\Gamma(\alpha,\,\beta)$ in $\Gamma(\alpha-2,\,\beta),\,\Gamma(\alpha-4,\,\beta),\,\ldots,\,\Gamma(\alpha_1,\,\beta),$ yields a graph $\Gamma_1\in\Gamma(\alpha_1,\,\beta).$ But the class of graphs in $\Gamma(\alpha_1,\,\beta)$ are isomorphic to the class of graphs in $\Gamma(\beta,\,\alpha_1),$ that is, to each graph of $\Gamma(\alpha_1,\,\beta)$ there is an isomorphic copy of it in $\Gamma(\beta,\,\alpha_1)$ and vice versa. Hence $\Gamma_1$ can be considered as a graph $\Gamma_2\in\Gamma(\beta,\,\alpha_1).$ Then successive reduced graphs of $\Gamma_2$ yields a graph $\Gamma_3\in\Gamma(\beta_1,\,\alpha_1).$ As the class of graphs in $\Gamma(\beta_1,\,\alpha_1)$ are isomorphic to the class of graphs in $\Gamma(\alpha_1,\,\beta_1),$ $\Gamma_3$ can be considered as a graph  $\Gamma_4\in\Gamma(\alpha_1,\,\beta_1).$ We shall call $\Gamma_3\in\Gamma(\beta_1,\,\alpha_1)$ or $\Gamma_4\in\Gamma(\alpha_1,\,\beta_1)$ (note that $\Gamma_3\cong \Gamma_4$) as a reduced graph of $\Gamma\in\Gamma(\alpha,\,\beta)$ according to the circumstances.

Similarly, by saying, for $\Gamma^\prime\in\Gamma(\nu,\,\delta)$ obtain its lifted graph $\Gamma_{k+\ell}\in \Gamma(\nu+2k,\,\delta+2\ell),$ where $k$ and $\ell$ are not simultaneously zero, we mean the following: for $\Gamma^\prime\in\Gamma(\nu,\,\delta),$ by lifting, we obtain a graph $\Gamma_1\in\Gamma(\nu+2,\,\delta).$ Similarly, for $\Gamma_1\in\Gamma(\nu+2,\,\delta),$ by lifting, we obtain $\Gamma_2\in\Gamma(\nu+4,\,\delta).$ Successively, we can get the graph $\Gamma_k\in\Gamma(\nu+2k,\,\delta).$  But the class of graphs in $\Gamma(\nu+2k,\,\delta)$ are isomorphic to the class of graphs in $\Gamma(\delta,\,\nu+2k),$ that is, to each graph of $\Gamma(\nu+2k,\,\delta)$ there is an isomorphic copy of it in $\Gamma(\delta,\,\nu+2k)$ and vice versa. Hence $\Gamma_k$ can be considered as a graph of $\Gamma(\delta,\,\nu+2k).$  From this, by successive liftings, we obtain $\Gamma_{k+1}\in \Gamma(\delta+2,\,\nu+2k),\,\Gamma_{k+2}\in\Gamma(\delta+4,\,\nu+2k),\,\ldots,\,\Gamma_{k+\ell}\in\Gamma(\delta+2\ell,\,\nu+2k).$ We call $\Gamma_{k+\ell}$ as a lifted graph of $\Gamma^\prime\in\Gamma(\nu,\,\delta).$ As the two classes of graphs $\Gamma(\nu+2k,\,\delta+2\ell)$ and $\Gamma(\delta+2\ell,\,\nu+2k)$ have the same set of graphs, upto isomorphism, we can consider $\Gamma_{k+\ell}$ as an element of $\Gamma(\nu+2k,\,\delta+2\ell).$


\end{rmk}

The idea in the next theorem is based on \cite{JCTB.46.142-153}.

\begin{thm}\label{THM.directcirculant}
Every $4$-regular connected circulant graph $\Gamma\in\Gamma(\alpha,\,\beta)$ of even order with jumps of different parity has the property $Q_2$ with respect to a Hamilton cycle decomposition.
\end{thm}
\NI{\bf Proof.} Let $\Gamma\in\Gamma(\alpha,\,\beta).$ As observed in Lemma \ref{LEM1.directcirculant}, $\alpha$ and $\beta$ are of different parity and hence we assume that $\alpha$ is even and $\beta$ is odd. Obtain the reduced graph $\Gamma^\prime\in\Gamma(1,\,2)$ of $\Gamma,$ see Remark \ref{DC.remark4}. If $\Gamma^\prime$ is simple, then it has a Hamilton cycle decomposition satisfying properties $Q_1$ and $Q_2,$ by Lemma \ref{LEM2.directcirculant}. Now we \lq\lq lift" $\Gamma^\prime$ to $\Gamma^{\prime\prime}\in\Gamma(\beta,\,2).$  By Lemma \ref{LEM1.directcirculant}, $\Gamma^{\prime\prime}$ has the properties $Q_1$ and $Q_2$ with respect to a Hamilton cycle decomposition of it. $\Gamma^{\prime\prime}$ can be considered as a graph $\Gamma^{\prime\prime\prime}\in\Gamma(2,\,\beta).$ Again lift $\Gamma^{\prime\prime\prime}$ to a graph in $\Gamma(\alpha,\,\beta),$ which is precisely $\Gamma,$ having the properties $Q_1$ and $Q_2$ with respect to a Hamilton cycle decomposition of it, by Lemma \ref{LEM1.directcirculant}. Thus the result is true if $\Gamma^\prime$ is simple.

Next we assume that $\Gamma^\prime$ is not simple. Let $\Gamma_1\in\Gamma(\alpha_1,\,\beta_1)$ be a reduced graph of $\Gamma\in\Gamma(\alpha,\,\beta)$ such that $\Gamma_1$ is simple but the reduced graph of $\Gamma_1$ in $\Gamma(\alpha_1-2,\,\beta_1)$ is not simple. The successive reduced graphs of $\Gamma_1$ results in $\Gamma^\prime\in\Gamma(2,\,1).$ Because of the above process of reductions and liftings of the classes of graphs, $\alpha_1$ and $\beta_1$ may have different parity from the parity of $\alpha$ and $\beta,$ respectively. $\Gamma^\prime$ is not simple only in three cases, see \cite{JCTB.46.142-153}:

(i) $\alpha_1=3$ and $c=0;$ in that case we obtain loops in $\Gamma^\prime.$

(ii) $\alpha_1=4$ and $c=0;$ in this case the edges in $\Gamma^\prime$ between $C_0^\prime$ and $C_1^\prime$ and, between $C_1^\prime$ and $C_0^\prime$ are the same, that is, they are multiple edges.

(iii) $\alpha_1=3$ and $c=\frac{k_a}{2};$ where $k_a$ is the order of the element $a$ in the generating set of $\Gamma_1;$ in this case we get multiple edges in $\Gamma^\prime.$

As the jumps of $\Gamma_1$ are of different parity, the first two cases correspond to the Cartesian product of cycles; in that case $c=0$ and $k_a=\beta_1.$  If $\beta_1\ge 5,$ again we reduce the graph $\Gamma_1\in\Gamma(3,\,\beta_1)$ or $\Gamma_1\in\Gamma(4,\,\beta_1)$ to a graph in $\Gamma(3,\,4)$ or $\Gamma(4,\,3).$ But $\Gamma(3,\,4)$ and $\Gamma(4,\,3)$ are isomorphic classes of graphs. Hence it is enough to consider $\Gamma(3,\,4)$ with $c=0.$ Then the class of graphs $\Gamma(3,\,4)$ with $c=0$ reduces to the single graph $C_3\Box C_4,$ where $\Box$ is the Cartesian product of graphs. In Figure 12, a Hamilton cycle decomposition is shown, where $C=(1\,2\,3\,4)$ is an odd alternating $4$-cycle.

\vspace{-.4cm}
\begin{center}
\includegraphics{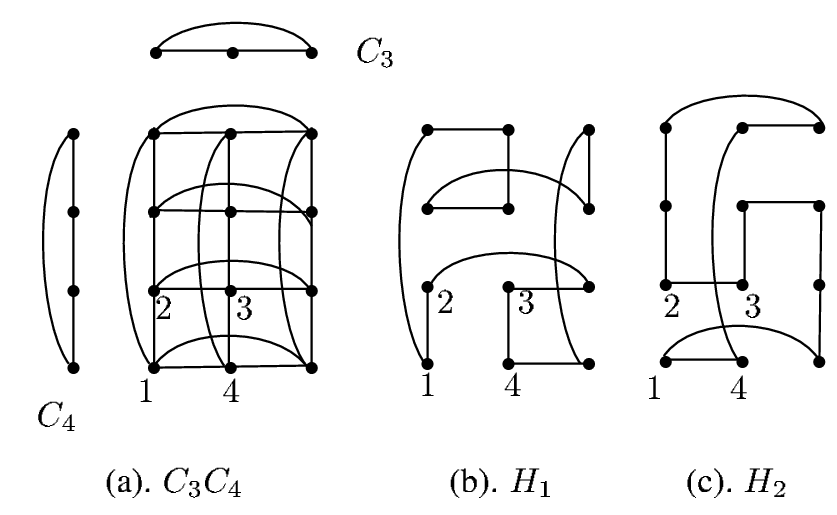}
\end{center}
\vspace{-.4cm}

{\small The graph $C_3\Box C_4\in\Gamma(3,\,4)$ and a Hamilton cycle decomposition $\{H_1,\,H_2\}$ with property $Q_2.$}\vspace{-.4cm}\begin{center}Figure 12
\end{center}

Next we consider the last case $\alpha_1=3,\,c=\frac{k_a}{2}.$ As $\Gamma_1$ is simple, then $k_a\ge 3$ and $\beta_1\ge 2,$ as $c=\beta_1=\frac{k_a}{2}\ge2.$ If $\beta_1\ge 4,$ we first reduce the graph $\Gamma_1$ to a graph in $\Gamma(3,\,2).$ In this case there is only one graph in $\Gamma(3,2),$ see \cite{JCTB.46.142-153},
 as shown in Figure 13. In the Hamilton cycle decomposition shown in Figure 13, $C=(1\,2\,3\,4)$ is an odd alternating $4$-cycle.

\begin{center}
\includegraphics{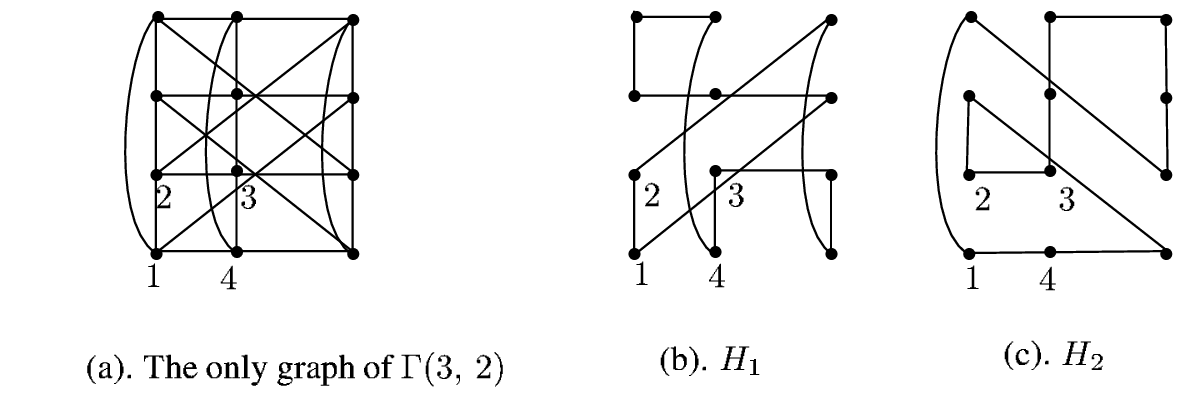}
\end{center}

\vspace{-.4cm}
{\small The only graph of $\Gamma(3,\,2)$ and a Hamilton cycle decomposition with property $Q_2.$}\vspace{-.4cm}\begin{center}Figure 13
\end{center}

\vspace{-.13cm}
To complete the proof, first we lift the graph $C_3\Box C_4\in\Gamma(3,\,4)$ of Figure 12 to the graph $G\in\Gamma(\beta,\,4)$ with properties $Q_1$ and $Q_2$ with respect to a Hamilton cycle decomposition, by Lemma \ref{LEM1.directcirculant}. Let $G_1\in\Gamma(4,\,\beta)$ be isomorphic to $G.$ As $G$ has the properties $Q_1$ and $Q_2,$ so does $G_1.$ Now we lift the graph $G_1\in\Gamma(4,\,\beta)$ to the graph $\Gamma\in\Gamma(\alpha,\,\beta).$ Clearly, $\Gamma\in\Gamma(\alpha,\,\beta)$ has a Hamilton cycle decomposition, by Lemma \ref{LEM1.directcirculant}. Similarly, we prove it for the graph in Figure 13.

This completes the proof of theorem.\hfill$\Box$

The following theorem of Jha \cite{Comput.Mat.Appl.31.11-19} is used in our proof of Theorem \ref{THM.CIR.1}.

\begin{thm}\emph{\cite{Comput.Mat.Appl.31.11-19}}
\label{THM.MAIN.directcirculant}
Let $\{H_1,\,H_2\}$ be a Hamilton cycle decomposition of a $4$-regular graph $G$ of even order $m$ containing an odd alternating four cycle, that is property $Q_2,$ with respect to $\{H_1,\,H_2\},$ then the graph $C_n\times G,\,n$ even, admits a Hamilton cycle decomposition.\hfill$\Box$
\end{thm}

\NI{\bf Proof of Theorem \ref{THM.CIR.1}.}

Let $|V(G)|=m$ and $|V(H)|=n.$ Since $H$ is Hamilton cycle decomposable, we have $H=C_n^1\oplus C_n^2\oplus \ldots\oplus C_n^\ell,$ where each $C_n^i$ is a Hamilton cycle of $H.$ If at least one of $G$ or $H$ has odd order, then $G\times H$ is Hamilton cycle decomposable, by Theorem \ref{Ann.DM.3.21-28}, and hence assume that both $G$ and $H$ are of even order. Since $G$ has the property $Q,$ $G$ can be decomposed into $4$-regular connected circulants, that is, $G=G_1\oplus G_2\oplus \ldots \oplus G_{\ell^\prime},$ where each $G_i$ is a $4$-regular circulant graph with jumps of different parity. Thus $G_i,\,1\le i\le \ell^\prime,$ can be decomposed into two Hamilton cycles with property $Q_2,$ by Theorem \ref{THM.directcirculant}. Now $G\times H\cong (G_1\oplus G_2\oplus \ldots\oplus G_{\ell^\prime})\times (C_n^1\oplus C_n^2\oplus \ldots\oplus C_n^\ell)=(G_1\times C_n^1)\oplus \ldots (G_1\times C_n^\ell)\oplus \ldots \oplus (G_{\ell^\prime}\times C_n^1)\oplus \ldots (G_{\ell^\prime}\times C_n^\ell).$ But each $G_i\times C_n^j,\,1\le i\le \ell^\prime,\,1\le j\le \ell,$ can be decomposed into Hamilton cycles, by Theorem \ref{THM.MAIN.directcirculant}.

This completes the proof of the theorem.\hfill$\Box$

%

\NI {\bf Conclusion.} In \cite{{DM.268.49-58}, {IJPAM.23.723-729}, {ARS.80.33-44},  {DM.308.3586-3606}, {DM.310.2776-2789}}, existence of Hamilton cycle decompositions of the graphs $K_r\times K_s,\,K_r\times K_{s,\,s},\,K_{r(s)}\times K_{m(n)},\,K_{r,\,r}\times K_{m(n)}$ are proved; the factor graphs in the product graphs are either complete or complete multipartite graphs, which are very dense graphs. However, if we consider one of the factor graphs as circulant, in $G\times H,$ with same number of odd and even jumps, irrespective of the number of jumps, which can be paired so that the resulting set of edges induces connected $4$-regular circulants, then $G\times H$ is Hamilton cycle decomposable. This proves that with this additional condition on $G$ or $H,$ the factor graphs $G$ and $H$ need not be dense. For example, Hamilton cycle decomposition of $(K_{4r+2}-F)\times H,$ where $F$ is a $1$-factor of $K_{4r+2},$ follows from our Theorem \ref{THM.CIR.1} whenever $H$ is Hamilton cycle decomposable multigraph. Also one can conclude that tensor products of certain sparse Hamilton cycle decomposable circulant graphs are Hamilton cycle decomposable.

\NI{\bf Acknowledgments:} The authors would like to thank the Department of Science and Technology, Government of India, New Delhi, for partial financial assistance through Grant No: SR/S4/MS:481/07.

\end{document}